\documentclass[11pt]{amsart}

\usepackage{fullpage, graphicx, enumitem}
\usepackage{url}
\usepackage{verbatim}
\usepackage{bbm}
\usepackage{amssymb,amsmath}
\usepackage{amsfonts,savesym,graphicx,bm,amsthm}
\usepackage{color}
\usepackage[mathscr]{eucal}
\usepackage[all]{xy}

\definecolor{hot}{RGB}{65,105,225}
\usepackage[pagebackref=true,colorlinks=true, linkcolor=hot ,  citecolor=hot, urlcolor=hot]{hyperref}

\newtheorem{theorem}{Theorem}[section]
\newtheorem{lemma}[theorem]{Lemma}

\newtheorem{theorem-definition}[theorem]{Theorem-Definition}
\newtheorem{corollary}[theorem]{Corollary}
\newtheorem{proposition}[theorem]{Proposition}
\newtheorem{prop}[theorem]{Proposition}

\newtheorem{definition-theorem}[theorem]{Definition-Theorem}

\newtheorem{theorem-defintion}[theorem]{Theorem-Definition}
\newtheorem{corollary-definition}[theorem]{Corollary-Definition}
\newtheorem{definition-proposition}[theorem]{Definition-Proposition}

\theoremstyle{definition}
\newtheorem{example}[theorem]{Example}
\newtheorem{definition}[theorem]{Definition}

\newtheorem{remark}[theorem]{Remark}

\newtheorem*{ack}{Acknowledgements}

\newtheorem{subs}[theorem]{ }

\numberwithin{equation}{section}

\def\bC{\mathbb{C}}
\def\be{\begin{equation}}
	\def\ee{\end{equation}}
\def\Exp{\mathrm{Exp}}
\def\bN{\mathbb{N}}
\def\bZ{\mathbb{Z}}
\def\al{\alpha}
\def\bQ{\mathbb{Q}}
\def\lam{\lambda}

\def\ord{\text{ord}}
\def\Spec{\text{Spec}}
\def\xa{\xrightarrow}
\def\pa{\partial}
\newcommand{\ubul}{{\,\begin{picture}(-1,1)(-1,-3)\circle*{2}\end{picture}\ }}
\def\cI{\mathcal{I}}
\def\cJ{\mathcal J}
\def\cL{\mathcal L}
\def\cV{\mathcal V}
\def\cA{\mathcal A}
\def\divi{\mathrm{div}}
\def\Supp{\mathrm{Supp}}
\def\bs{{\bm s}}

\author{Nero Budur}
\address{Department of Mathematics, KU Leuven, Celestijnenlaan 200B, 3001 Leuven, Belgium;  YMSC, Tsinghua University, 100084 Beijing, China;  BCAM, Mazarredo 14, 48009 Bilbao, Spain.}
\email{nero.budur@kuleuven.be}

\author{Quan Shi}
\address{Department of Mathematical Sciences, Tsinghua University, Beijing, 100084, P. R. China.}
\email{shiq24@mails.tsinghua.edu.cn / thusq20@gmail.com}

\author{Huaiqing Zuo}
\address{Department of Mathematical Sciences, Tsinghua University, Beijing, 100084, P. R. China.}
\email{hqzuo@mail.tsinghua.edu.cn}

\title{Polar loci of multivariable archimedean zeta functions}

\begin{document}

	\begin{abstract}
		We determine, up to exponentiating, the polar locus of the multivariable archimedean zeta function associated to a finite collection of polynomials $F$. The result is the monodromy support locus of $F$, a topological invariant. We give a relation between the multiplicities of the irreducible components of the monodromy support locus and the polar orders. These generalize results of Barlet for the case when $F$ is a single polynomial. Our result determines the slopes of the polar locus of the zeta function of $F$,  closing a circle of results of Loeser,   Maisonobe, Sabbah. We apply our main result to elucidate the topological information contained by the oblique part of the zero locus of any ideal of Bernstein-Sato type.
		
	\end{abstract}

	\maketitle

	\tableofcontents

	\section{Introduction} \label{Intro}

	Let $X$ be a complex manifold of dimension $n$ and let $F = (f_1,...,f_r)$ be an ordered collection of holomorphic functions on $X$ such that $f:=\prod_{i=1}^rf_r\neq 0$  has non-empty zero locus.  For example this can be a smooth irreducible complex algebraic variety and $f_i$ polynomials, or,  $F$ can be the representative of a germ of a complex analytic mapping on a complex manifold. These two are particular cases of the subclass defined by the condition:
	\be\label{eqCon}
	\tag{*}\text{The zero locus of }f=\prod_{i=1}^rf_i\neq 0\text{ is non-empty and admits a finite Whitney stratification.}
	\ee 	
	We shall assume throughout the paper this condition is satisfied.

	The \textit{multivariable archimedean zeta function of $F$ and $\omega$}, where $\omega$ is compactly supported smooth $(n,n)$-form on $X$, is  the meromorphic extension to $\mathbb C^r$ of 
	\begin{displaymath}
		Z_{F,\omega}(s_1,...,s_r)=\int_{X} \vert f_1\vert^{2s_1} \cdots \vert f_r\vert^{2s_r} \omega,
	\end{displaymath}
	which is a priori defined for $\mathrm{Re}(s_i)> 0$ for all $i$. By the {\it multivariable archimedean zeta function of $F$} one means the distribution $Z_F$ defined by $Z_{F,\omega}$ as $\omega$ varies. 
	
	The case $r=1$ has been extensively studied, see  \cite{LC22} for a survey. The case $r>1$ has probably been first studied in connection with Feynman amplitudes, see the introduction of \cite{BWZ}. Feynman integrals are integrals that fall within the type $Z_{F,\omega}$, where the choice of compactly-supported $\omega$  corresponds to a cut-off regularization. There is  an extensive literature on Feynman integrals, see the book \cite{We}. From a physicist's point of view, the goal with these integrals is to compute them when defined and understand their asymptotic behaviour as $s$ approaches the poles of $Z_{F,\omega}(s)$. A first step is thus to determine the polar locus. This is the subject of this article.

	The {\it polar locus of $Z_F$} is defined to be the union over all $\omega$ of the polar loci of $Z_{F,\omega}$.  The polar locus of $Z_F$ is contained in a union of hyperplanes  in $\bC^r$ of the form $\sum_{i=1}^r c_i s_i + c_0 = 0$, where all $c_i \in \mathbb N$ and $c_0 > 0$. This can be seen by two different methods, as in the case $r=1$: the existence of a log resolution of $f$, and the existence of solutions to functional equations of Bernstein-Sato type, see \cite{Sab87, Loe89, Gy}.

	The determination of the polar locus of $Z_F$ is an open problem. Set
	$$\mathrm{Exp} : \mathbb C^r \to (\mathbb C^*)^r,\quad (\al_1,...,\al_r) \mapsto (e^{2\pi i \al_1},...,e^{2\pi i\al_r}).$$
	Define the {\it slope} of the hyperplane  $\sum_{i=1}^r c_i s_i + c_0 = 0$ to be $[c_1:...:c_r] \in \mathbb P^{r-1}$. 
	Upper bounds on the polar locus, its  image under $\Exp$, and on the set of its slopes where given in \cite{Sab87, Loe89, Mai}. The case $r=2$ was studied under further conditions in \cite{BM}.
	
	In this paper we  determine the image under $\Exp$ and the set of slopes of  the polar locus of $Z_{F}$. 
	
	To state the results, we recall that one can attach to $F$ two related invariants, one topological and the other analytical. The first is the {\it monodromy support} $S(F)\subset (\bC^*)^r$, see \ref{Alexander module and its specialization}. In the case $r=1$, that is, when $F$ is a single holomorphic function $f$, the monodromy support $S(f)$ is the set of all eigenvalues of the monodromy on the cohomologies of the Milnor fibers of $f$. For $r>1$, there is locally a monodromy around each divisor $\{f_i=0\}$, and $S(F)$ captures these actions simultaneously. The monodromy support $S(F)$ can be computed from log resolutions of $f$ by  \cite[Theorem 1.3]{BLSW17}, which generalizes the formula of A'Campo \cite{AC} for computing  $S(f)$. For example, if each $f_j$ is a product of linear forms, $S(F)$ admits an explicit combinatorial formula, see \cite[Theorem 1.2]{BLSW17}, \cite[6.7, 6.8]{Bud15}.

	A second invariant is the \textit{Bernstein-Sato ideal} $B_F\subset\bC[s_1,\ldots,s_r]$, see \ref{B-S polynomial}. This is a non-zero ideal first defined locally at every point of the zero locus of $f$. Then the assumption (\ref{eqCon}) is used to define it globally as the intersection of the local Bernstein-Sato ideals.  In the case $r=1$, $B_f$ is principal and generated by the classical Bernstein-Sato polynomial of $f$, denoted by $b_f(s)$. In the algebraic case, $B_F$ can be in principle computed with a computer algebra system. We denote by $Z(B_F)\subset \bC^r$ the zero locus of the Bernstein-Sato ideal.
	
	The structure of these two invariants and their relationship is given by:

	\begin{theorem}\label{Zero_Loci_of_B-S}
		Let $X$ be a complex manifold   and $F = (f_1,...,f_r) : X \to \mathbb C^r$ be an analytic mapping satisfying (\ref{eqCon}). Then:
		
		\begin{enumerate}
			
			\item (Budur-Liu-Saumell-Wang \cite{BLSW17}) $S(F)$ is a conjugation-invariant  finite-union of torsion-translated complex affine subtori of codimension one in $(\mathbb C^*)^r$.
			\item (Sabbah \cite{Sab87}, Gyoja \cite{Gy}) $B_F$ contains a product of polynomials of the  type $\sum_{i=1}^rc_is_i+c_0=0$ with all $c_i\in\bN$ and $c_0>0$. In particular, every irreducible component of $Z(B_F)$ is contained in a  hyperplane of this type.
			\item (Maisonobe \cite{Mai}) Every irreducible component of $Z(B_F)$ of codimension $>1$ can be translated by an element of $\bZ^r$ into a component of codimension one.
			\item 
			(Budur-van der Veer-Wu-Zhou \cite{Zero_Loci_of_Bernstein-Sato_Ideals}) $
			\mathrm{Exp}(Z(B_F)) = S(F).
			$		
		\end{enumerate}
	\end{theorem} 
	Here, a complex affine subtorus of $(\bC^*)^r$ is an algebraic subgroup $(\bC^*)^p\subset (\bC^*)^r$ for some $p\le r$. Part (4) will be crucial for this paper, another proof of it is in \cite{VV}.
	
	When $r=1$,  Theorem \ref{Zero_Loci_of_B-S} is classical, due to Malgrange, Kashiwara \cite{Rationality_of_Roots_of_B-Function_Kashiwara, Kashiwara_Malgrange_Theorem_Malgrange, Kashiwara_Malgrange_Theorem_Kashiwara}.  
	
	We recall the following collection of results on the archimedean zeta functions $Z_f$ in this case:
	
	\begin{theorem}\label{Triangle of three local invariant}
		Let $f : X \to \mathbb C$ be a holomorphic function on a complex manifold satisfying (\ref{eqCon}). Then:
		
		\begin{enumerate}
			
			\item (Bernstein \cite{Be}, Bj\"ork \cite{Bj},  Kashiwara \cite{Rationality_of_Roots_of_B-Function_Kashiwara}, Igusa \cite{Igu}) If  $\alpha$ is a pole of $Z_{f}(s)$ then  $\alpha \in \mathbb Q_{<0}$ and it is a root of $\prod_{k = 0}^{\lceil-\alpha\rceil-1} b_f(s+k)$. Moreover, the order of $\alpha$ as a pole can be at most the multiplicity as a root.
			
			
			\item (Barlet \cite{Bar84}) If $\lambda$ is a monodromy eigenvalue for $f$, that is, if $\lambda \in S(f)$, then there exists   a pole $\al$ of  $Z_{f}(s)$  such that $\lambda = e^{2\pi i \al}$.  In particular,		\begin{displaymath}
				\{ \lambda\in\bC^*\mid  \lambda=e^{2\pi i \al}\text{ for some pole }\al\text{ of }Z_{f}\} = S(f).
			\end{displaymath}
			
			\item  (Barlet) Suppose $\lam\in S(f)$ is an eigenvalue of the monodromy on the cohomology $H^j$ of the Milnor fiber of $f$ 	
			at a point $x \in \{f = 0\}$, with a Jordan block of size  $m\ge 1$. Then:
			\begin{enumerate}
				\item (\cite{Bar84}) 
				There exists  a pole $\al$ of $Z_{f}$ of order $\geq m$ such that $e^{2\pi i\al} = \lam$. 
				\item (\cite{Bar86}) Write $\lam=e^{2\pi i\al}$ with $\al\in (-1,0]\cap \bQ$. If $j\ge 1$ then $\al-j$ is a pole of $Z_f$ of order $\ge m$. 
				
				\item (\cite{Bar84b}) If $\lam=1$ and $j\ge 1$, there exists  a negative integer pole  of $Z_{f}$ of order $\geq m+1$.
			\end{enumerate}
		\end{enumerate}
	\end{theorem}
	
	The above results are originally phrased locally, but  it is easy to globalize them, see Section \ref{Preliminary}.
	
	Barlet's results determine thus completely the image under $\Exp$ of poles of the archimedean zeta function $Z_{f}$, and this image is a topological invariant. 
	In contrast, it was shown recently  that not every root $\al\in Z(B_f)$ is a pole of the archimedean zeta function $Z_{f}$ by Davis,  L\"orincz, and Yang \cite{DLY24}. 
	
	Theorem \ref{Triangle of three local invariant} is one of the motivations behind the Monodromy Conjecture of Igusa \cite{Ig}, Denef-Loeser \cite{DL98} for  $p$-adic and, respectively, motivic analogs of archimedean zeta functions.

	\subs{\bf Main results.}
	In this paper we  generalize Theorem \ref{Triangle of three local invariant} to the case $r>1$. Part (1) of Theorem \ref{Triangle of three local invariant} generalizes easily, see Proposition \ref{B-S proof of one direction}. In particular, it gives $\Exp(P(Z_F))\subset \Exp(Z(B_F)) = S(F)$, where $P(Z_F)$ is the polar locus of the multivariable archimedean zeta function $Z_{F}$. The generalization of part (2)  of Theorem \ref{Triangle of three local invariant} is the reverse inclusion:
	
	\begin{theorem}\label{Thm_Support_Pole}
		Let $X$ be a complex manifold   and $F = (f_1,...,f_r) : X \to \mathbb C^r$ be an analytic mapping satisfying (\ref{eqCon}). Then 	in $(\bC^*)^r$,	we have	$
		\mathrm{Exp}(P(Z_{F})) = S(F).
		$
	\end{theorem}
	
	In particular, $\Exp(P(Z_F)$ is a topological invariant, it can also be computed from log resolutions of $f$, and, in case $f$ is a product of linear polynomials, it is a combinatorial invariant. This also gives a new proof of the codimension-one claim in Theorem \ref{Zero_Loci_of_B-S} (1), due to  \cite{BLSW17}.

	Recall that we have the notion of slope for the polar hyperplanes of $Z_F$. This notion is trivial if $r=1$, but it is non-trivial if $r>1$.
	Theorem \ref{Thm_Support_Pole} has a consequence that the slopes of the polar locus of $Z_F$ are also completely determined. In fact one has several interpretations, topological, algebraic, and geometric:
	
	\begin{theorem}\label{Slopes of zeta function}
		Let $X$ be a complex manifold   and $F = (f_1,...,f_r) : X \to \mathbb C^r$ be an analytic mapping satisfying (\ref{eqCon}).  Then the following sets are equal:
		\begin{enumerate}[label=(\alph*)]		
			\item the  slopes of the polar locus of $Z_F$,
			\item the  slopes of the tangent hyperplanes at the identity of the complex affine subtori of $(\bC^*)^r$ that are translates of the irreducible components of the monodromy support locus $S(F)$,
			\item the  slopes of the irreducible components of the zero locus $Z(B_F)$ of the Bernstein-Sato ideal,
			\item the  slopes of the irreducible components of the zero locus of the initial ideal of $B_F$ (generated by the highest-degree terms of polynomials in $B_F$),
			\item the  slopes of the irreducible components of $\pi_2(W_{F,0}^\sharp)$, where:
			\begin{itemize}
				\item $p : T^*X \to X$ is the cotangent bundle of $X$,
				\item  $W_F^\sharp \subset T^*X \times \mathbb C^r$ is the closure of	$$
				\left\{\left(x,\textstyle\sum_{i=1}^r\alpha_i({df_i}/{f_i})(x), \alpha_1,\ldots, \alpha_r\right)\in T^*X\times \mathbb C^r\mid f(x)\ne 0\right\},
				$$
				\item $W^{\sharp}_{F,0} := W_F^{\sharp} \cap p^{-1}(\{f = 0\})$, and $\pi_2$ is the projection to $\bC^r$.
			\end{itemize}
		\end{enumerate}
	\end{theorem}

	Here, $\pi_2(W^{\sharp}_{F,0})$ is a union of central hyperplanes by \cite{BMM}, the set $W^{\sharp}_{F,0}$ having been introduced earlier by Kashiwara and Kawai as the characteristic variety of a certain $\mathscr D$-module. The equivalence between (c), (d), (e) is due to Maisonobe \cite{Mai}. The equivalence between (b) and (c)-(e) follows from Theorem \ref{Zero_Loci_of_B-S}  above. So only the equivalence between (a) and (b) is new here, and it follows from Theorem \ref{Thm_Support_Pole}.

	Theorem \ref{Slopes of zeta function} closes a  circle of results on upper bounds for the slopes of the polar locus of $Z_F$ due to \cite{Sab87,Loe89,Mai}. It says in fact that Maisonobe's upper bound is exactly the set of slopes of the polar locus of $Z_F$. Using the theorem we correct one, and improve both, of the two main examples from \cite{BM}, see Examples \ref{exaBM}, \ref{exaBM2}.	
	
	We will give yet another characterization of the oblique slopes of the polar locus, see \ref{subObl} below.

	The two sides of the equality in Theorem \ref{Thm_Support_Pole} are not only set-theoretically related.  On one hand, for a polar hyperplane of the archimedean zeta function $Z_F$ we have the concept of polar order. On the other hand, a problem we encountered is to define a good notion of order for the components of the monodromy support locus. By definition, $S(F)$ is the union  of the monodromy support loci $\mathrm{Supp}_x^i(\psi_F\mathbb C_X)$ of the cohomology sheaves $\mathcal H^i$ of stalks at $x\in\{f=0\}$ of Sabbah's specialization complex $\psi_F\bC_X$, see  \ref{Alexander module and its specialization}. This is a complex of sheaves of modules over $A=\bC[t_1^{\pm},\ldots,t_r^{\pm}]$, which for  $r=1$ agrees with the shift  $[-1]$ of Deligne's nearby cycles complex  endowed with the monodromy action.
	To an irreducible component $C$ of codimension one  of $\mathrm{Supp}_x^i(\psi_F\mathbb C_X)$  we associate a multiplicity $\mathrm{ord}_{C}^i(F,x)$ as the order of vanishing of $C$ in the annihilator in $A$ of $\mathcal H^i((\psi_F\bC_X)_x)$. In case $r=1$, this  notion coincides with the maximal size of Jordan blocks of with prescribed eigenvalue of the monodromy on the cohomology $H^{i-1}$ of the Milnor fiber of $f$ at $x$.   
	
	Our generalization of Theorem \ref{Triangle of three local invariant} (3) is given in Theorem \ref{thmGenbc}. 	We state here only  the least technical part of it, which is
	our generalization of Theorem \ref{Triangle of three local invariant} (3) (a):

	\begin{theorem}\label{Thm_Order controlling}
		Let $X$ be a complex manifold of dimension $n$  and $F = (f_1,...,f_r) : X \to \mathbb C^r$ be an analytic mapping satisfying (\ref{eqCon}). 
		Let $C$ be an irreducible component of codimension one of $\Supp^i_x(\psi_F\bC_X)$. Then $C$ is the $\Exp$ image of  a polar hyperplane of $Z_F$ of polar order $\ge \mathrm{ord}_{C}^i(F,x)$.
	\end{theorem}
	
	Together with Proposition \ref{B-S proof of one direction} (3) and {Lemma \ref{propagation_of_poles}}, this implies a relation between the orders of vanishing along codimension one components of $B_F$ and $S(F)$:

	\begin{corollary}
		With the assumptions as in Theorem \ref{Thm_Order controlling}, if $H \subset Z(B_F)$ is a hyperplane such that $\mathrm{Exp}(H) = C$, then $\sum_{k \in \mathbb Z} \mathrm{ord}_{H-k\cdot \bm 1}(B_F) \geq \mathrm{ord}_{C}^i(F,x)$.
	\end{corollary}

	To prove our results we use two  techniques. Firstly,  polar hyperplanes are compatible with the specialization of $F=(f_1,\ldots,f_r)$  to a product of generic powers of $f_1,\ldots ,f_r$ in the following sense:

	\begin{theorem}\label{Specialization of Archimedean Zeta Function}
		Let $X$ be a complex manifold of dimension $n$  and $F = (f_1,...,f_r) : X \to \mathbb C^r$ be an analytic mapping satisfying (\ref{eqCon}). 
		\begin{enumerate}	
			\item Suppose  $\sum_{i=1}^r c_i s_i + c_0 = 0$ with $c_0,...,c_r\in \mathbb N$, is a polar hyperplane of $Z_{F}(s_1,...,s_r)$ of order $m$. Then there exists an analytic open dense  subset $\Omega$ of $\sum_{i=1}^r s_i = 1$ such that $-c_0/\sum_{i=1}^r c_i b_i$ is a pole of order $m$ of $Z_{f_1^{b_1}\cdots f_r^{b_r}}(s)$ for all $(b_1,...,b_r) \in \mathbb N_{>0}^r$ with $(b_1,\ldots, b_r)/{\sum_{i=1}^r b_i} \in \Omega$.
			
			\item  Let $(b_1,...,b_r) \in \mathbb N_{>0}^r$. If $\alpha$ is a pole of $Z_{f_1^{b_1}\cdots f_r^{b_r}}(s)$, then there exists a polar hyperplane $\sum_{i=1}^r c_i s_i + c_0 = 0$ of $Z_{F}(s_1,\ldots,s_r)$ such that		 $\sum_{i=1}^r c_ib_i \alpha + c_0 = 0$. More generally, if $\alpha$ is a pole of order $m$ of $Z_{f_1^{b_1}\cdots f_r^{b_r}}(s)$, then the sum of the orders of all polar hyperplanes of $Z_{F}$ containing $(b_1\alpha,...,b_r\alpha)$ is $\geq m$.
		\end{enumerate}	
	\end{theorem}
	
	As an aside,  this together with  \cite[Theorem 1.7]{SZ25} give another proof  of the multivariable $n/d$-conjecture for complete splittings of hyperplane arrangements \cite{Wu22}, see Theorem \ref{thmNDConj}.

	There is a similar compatibility between forming the monodromy support locus and specializations of $F$ to the product of generic powers of $f_1,\ldots, f_r$, due to Sabbah \cite{Alexander Module}. Then the statement of 	Theorem \ref{Thm_Support_Pole} is reduced to the case $r=1$.

	Secondly, we arrived in a round-about way at the definition of $\ord_{C}^i(F,x)$. What we were looking for was a definition compatible with specializations as above. Since annihilators are not in general compatible with base change, no good notion of scheme structure on $\Supp^i_x(\psi_F\bC_X)$ was available. This led us to	introduce the \textit{cohomological determinantal-factor ideals} in Section \ref{Cohomological determinantal-factor ideal}. 
	These  are generalizations of Fitting ideals to complexes of modules and have interesting geometric interpretations.
	They refine the  cohomology jump ideals, see \cite{BW15jump}, but  the latter  do not provide the right scheme structure for our purposes either. 
	For comparison, the cohomology jump ideal $\mathcal J^i_k$ provides a closed subscheme structure on the locus of points where the $i$-degree cohomology has dimension $\ge k$. The cohomological determinantal-factor ideal $\mathcal I^i_k$ provides a new closed subscheme structure on the locus of points where the alternating sum of the dimensions of cohomology in degrees $\ge i$ is $\ge k+1$. In the context of a vector space $V$ with an endomorphism $t$, the ideals $\mathcal I^i_k$ of the associated $\bC[t]$-module $V$ recover the determinantal factors whose successive quotients define the invariant factors in the Smith normal form of $t$. 	
	
	We apply these new ideals to  $\psi_F\bC_X$. Denote  by $Z(\mathcal I_{k}^i(F,x))$ the zero locus in $(\bC^*)^r$ of the cohomological determinantal-factor ideal $\mathcal I_{k}^i(F,x) := \mathcal I_k^i((\psi_F\mathbb C_X)_x) \subset A = \mathbb C[t_1^{\pm},...,t_r^{\pm}]$.  Using the difference of divisors associated to $\mathcal I_{0}^i(F,x)$ and $\mathcal I_{1}^i(F,x)$, we define the multiplicities $\mathrm{ord}_{C}^i(F,x)$ in Definition  \ref{Application to Sabbah}. They turn out to coincide with the above definition in terms of the annihilator of $\mathcal H^i((\psi_F\bC_X)_x)$, see Proposition \ref{propAmu}.

	\begin{theorem}\label{zero locus of cohomological determinantal-factor ideals}\label{specialization of multiplicitiy}
		With the setup as in Theorem \ref{Thm_Support_Pole}, let $x\in \{f = 0\}$. 
		\begin{enumerate}
			\item		For $k,i\in \mathbb Z$, every irreducible component of $Z(\mathcal I_{k}^i(F,x))$ is a torsion-translated complex affine subtorus of $(\mathbb C^*)^r$.
			\item Let $C $ be an irreducible component of $\mathrm{Supp}_x(\psi_F\mathbb C_X)$ of codimension one. There exists a Zariski open dense subset $U$ of $C$ such that: (a) $U$ does not meet any other irreducible component of $S(F)$; (b) for all $(\lambda^{b_1},...,\lambda^{b_r})\in U$, where $\lambda\in\bC^*$ is a unit root and $b_1,...,b_r \in \mathbb N_{>0}$, the maximal size of Jordan blocks with eigenvalue $\lambda$ of the monodromy on the $(i-1)$-th cohomology of the Milnor fiber of $f_1^{b_1}\cdots f_r^{b_r}$ at $x$ equals $\mathrm{ord}_{C}^i(F,x)$.
%
%
%
%
		\end{enumerate}
		
	\end{theorem}
	
	With this the proof of Theorem \ref{Thm_Order controlling} is straight-forward.	
	
	To prove part (1) of Theorem \ref{zero locus of cohomological determinantal-factor ideals}, we relate  
	the subset $Z(\mathcal I_{k}^i(F,x))$  with cohomology jump loci of rank one local systems for which one disposes already of the torsion-translated subtori property.
	To state this, we need first a few general  properties of the cohomological determinantal-factor ideals for local systems. 
	
	Let $Y$ be a topological space with homotopy type of a finite CW complex. Denote by $\mathcal M_B(Y,m)$ the Betti moduli space of rank $m$ semi-simple $\bC$-local systems on $Y$. Define the \textit{cohomological determinantal-factor jump loci of $\mathcal M_B(Y,m)$} set-theoretically to  be
	\begin{displaymath}
		{\mathcal A}_k^i(Y,m) := \{\mathcal L \in \mathcal M_B(Y,m) \mid \sum_{l\geq i} (-1)^{l-i} \dim_\mathbb C H^l(Y,\mathcal L) \geq k\}.
	\end{displaymath}
	They are clearly determined by the boolean algebra generated by the {\it cohomology jump loci}
	$$
	{\mathcal V}_k^i(Y,m) := \{\mathcal L \in \mathcal M_B(Y,m) \mid  \dim_\mathbb C H^i(Y,\mathcal L) \geq k\}.
	$$
	
	\begin{proposition}\label{rigidity theorem of cohomological determinantal-factor jump loci}\label{rigidity of Aki}
		With notations are as above: 
		\begin{enumerate}
			\item
			For every $i,k \in \mathbb Z$, ${\mathcal A}_k^i(Y,m)$ is Zariski closed in $\mathcal M_B(Y,m)$. 
			\item If $Y$ is a smooth complex algebraic variety or the small ball complement of the germ of a complex analytic set, then ${\mathcal A}_k^i(Y,1)$ is a conjugation-invariant finite-union of torsion-translated complex affine subtori of the space of rank one local systems $\mathcal M_B(Y,1)$.
		\end{enumerate}
	\end{proposition}
	
	The cohomological determinantal-factor ideals place a closed subscheme structure on ${\mathcal A}_k^i(Y,m)$ which in general differs from the  structure induced from $
	{\mathcal V}_k^i(Y,m)
	$,
	cf. \cite{BW15jump}.

	Back to the situation of Theorem \ref{zero locus of cohomological determinantal-factor ideals}, we focus now on the rank one case $m=1$ and we drop it from the notation. Let $\nu_x:\mathcal M_B((\bC^*)^r)\to \mathcal M_B(U_x)$ be the pullback on local systems of rank one under the mapping $U_x\to (\bC^*)^r$ obtained by restriction from $F$, where  $U_x$ is the complement of $\{f=0\}$ in a small ball at $x$ in $X$. The local monodromy support loci   $\mathrm{Supp}_x(\psi_F\mathbb C_X)$ admit an interpretation
	in terms of cohomology support loci of rank one local systems on $U_x$, see Theorem \ref{torsion translated}. We refine this:

	\begin{proposition}\label{propRefW}
		With the setup as in Theorem \ref{zero locus of cohomological determinantal-factor ideals}:  	\begin{enumerate}
			\item
			The zero locus of the cohomology jump ideal $ \mathcal J_k^i((\psi_F\mathbb C_X)_x)$ is $\nu_x^{-1}(\mathcal V^i_{k}(U_x))$.
			
			\item  
			The zero locus of the cohomological determinantal-factor ideal $\mathcal I_k^i((\psi_F\mathbb C_X)_x)$ is $Z(\mathcal I_{k}^i(F,x))=\nu_x^{-1}(\mathcal A^i_{k+1}(U_x))$.
			
			\item $Z(\mathcal I_{0}^i(F,x))$ and $\Supp^i_x(\psi_F\bC_X)$ have the same codimension one irreducible components.
		\end{enumerate}
	\end{proposition}

	The monodromy support $\Supp^i_x(\psi_F\bC_X)$ is further related by a result of Liu and Maxim \cite[Proposition 3.7]{LM} with the local cohomology Alexander module of $F$.
	
	Proposition \ref{propRefW} gives the following for $r=1$:
	
	\begin{prop}\label{rmkAmuCon} Let $f:(\bC^n,0)\to (\bC,0)$ be the germ of a holomorphic function. Let $U$ be the complement of $f=0$ in a small ball.
		For $\lam\in\bC^*$, let $\cL_\lam$ be the rank one local system  on $U$ with monodromy $\lam$ around $f=0$. Let $F_{f,0}$ be the Milnor fiber of $f$ at the origin.
		Then for $i\in\bZ$:
		$$
		\{\text{eigenvalues of monodromy on }H^{i-1}(F_{f,x},\bC)\} =
		\{\lam\mid h^i(U_x,\cL_\lam)-h^{i+1}(U_x,\cL_\lam)+\ldots \ge 1\}.
		$$
	\end{prop}

	Recently Davis,  L\"orincz, and Yang \cite{DLY24} gave further refinements of Theorem \ref{Triangle of three local invariant} for $Z_f(s)$ in terms of Hodge filtrations. A natural next step would be to extend their  results to $Z_F(s_1,\ldots,s_r)$.
	
	\subs{\bf Application:  oblique slopes.}\label{subObl} The Bernstein-Sato ideal $B_F$ is one of the many ideals of Bernstein-Sato type $B_F^M$, which depend on
	a matrix $M\in \bN^{p\times r}$ with $p\ge 1$, see Definition \ref{defBFM}. One has $B_F=B_F^{\bm 1}$, where $\bm 1=(1,\ldots,1)\in\bN^r$. The ideals $B_F^M$ are non-zero. Also, $B_F\neq (1)$ if and only if for each row $1\le k\le p$ of $M$, there exists $1\le j\le r$ such that $f_j^{m_{kj}}$ is not invertible. The zero locus $Z(B_F^M)$ is contained in $Z(B_F)-\bN^r$, see Lemma \ref{lemPrep}. 
	
	The topological interpretation of $\Exp(Z(B_F))$ was generalized in \cite{BS2} to $\Exp(B_F^{\bm m})$ for $\bm m\in\bN^r$, see Section \ref{secObl}. What  $\Exp(Z(B_F^M))$ is for matrices $M$ with more than one row, is still open.  
	
	As an application of Theorem \ref{Thm_Support_Pole} we give a topological interpretation of the $\Exp$ image of the codimension one part of $Z(B_F^I)$, where $I\in\bN^{r\times r}$ is the identity matrix, of the  oblique part of $Z(B_F^M)$ for any $M$ with no zero rows,  and connect it with the  oblique part of the polar locus of $Z_F$.

	\begin{definition} The  hyperplane $\{c_1s_1+\ldots+c_rs_r+c_0=0\}$ in $\bC^r$, and its slope, are 
		called {\it  oblique} if $c_j\neq 0$ for all $1\le j\le r$.
		We denote by  $Z(B_F^{M})_{ob}$, $P(Z_F)_{ob}$, and $S(F)_{ob}$, the union of the 
		codimension one components with  oblique slopes of the zero locus $Z(B_F^{M})$, the polar locus $P(Z_F)$, and the monodromy support $S(F)$, respectively.
	\end{definition}
	
	\begin{theorem}\label{thmOBLIQUE}$\;$
		\begin{enumerate}
			\item All hyperplanes in $Z(B_F^I)$ are oblique. So, $Z(B_F^I)_{ob}$ is the codimension one part of $Z(B_F^I)$.
			\item	 For  all $M\in\bN^{p\times r}$ with no zero rows,
			$$
			\mathrm{Exp}(Z(B_F^{M})_{ob}) = \Exp(Z(B_F^I)_{ob})= \Exp(P(Z_F)_{ob})=\Exp(Z(B_F)_{ob})=S(F)_{ob}.
			$$ 
			\item Let $H=\{c_1s_1+\ldots+c_rs_r+c_0=0\}$ be a non-oblique hyperplane in $Z(B_F)$.
			Let $J$ be the ordered subset of $\{1,\ldots,r\}$ be such that $c_j\neq 0$ for all $j\in J$ and $c_j=0$ for $j\in\{1,\ldots,j\}\setminus J$. Then $H$ defines an oblique hyperplane of $Z(B_{F_J})$, where $F_J$ is the ordered tuple $(f_j\mid j\in J)$.
		\end{enumerate}
	\end{theorem}	
	
	\subs{\bf Layout.} Section \ref{Preliminary} is for preliminary material on monodromy support loci and Bernstein-Sato ideals $B_F$. In Section \ref{Cohomological determinantal-factor ideal} we introduce the cohomological determinantal-factor ideals. In Section \ref{secLocos} we apply them to prove Theorem \ref{specialization of multiplicitiy} and Propositions \ref{rigidity of Aki},  \ref{propRefW}, \ref{rmkAmuCon}. In Section \ref{pf1} we prove Theorems \ref{Thm_Support_Pole},   \ref{Thm_Order controlling},   \ref{Specialization of Archimedean Zeta Function}, \ref{thmGenbc}. In Section \ref{secObl}, we review the ideals $B_F^M$ of Bernstein-Sato type and prove Theorem \ref{thmOBLIQUE}. In Section \ref{secExa} we give  examples.

	\begin{ack}
		N. Budur was supported by a Methusalem grant and the grant G0B3123N from FWO.	Q. Shi and H. Zuo were supported by BJNFS Grant 1252009.  H. Zuo was supported by NSFC Grant 12271280. We thank Xiaoyu Su for discussions.
	\end{ack}

	\section{Preliminaries}\label{Preliminary}
	
	\subs{\bf Specialization complex.}\label{Alexander module and its specialization} 	
	Let $F = (f_1,...,f_r) : X \to \mathbb C^r$ be an analytic mapping from a complex manifold $X$ of dimension $n$ satisfying (\ref{eqCon}). Let $f = \prod_{i=1}^r f_i$ and $U = X\setminus f^{-1}(0)$. Let $j : U \hookrightarrow X$ and $i: f^{-1}(0) \hookrightarrow X$ be the natural embeddings. Consider the following commutative diagram of fiber products:
	\begin{displaymath}
		\xymatrix{
			f^{-1}(0) \ar@{^(->}^i[r] & X \ar_F[d] & \ar@{_(->}_j[l] U \ar_F[d] & \ar_\pi[l] \widetilde U \ar^{\widetilde F}[d]\\
			& \mathbb C^r & (\mathbb C^*)^r \ar@{_(->}[l] & (\widetilde{\mathbb C^*})^r \ar_{\widetilde\pi}[l]
		}
	\end{displaymath}
	where $\widetilde \pi = \mathrm{Exp} : (\widetilde{\mathbb C^*})^r=\bC^r \to (\mathbb C^*)^r$ is the universal covering map.
	\begin{definition}\label{defSpC}
		There is a well-defined functor 
		\begin{displaymath}
			\psi_F := i^{-1}Rj_*R\pi_!(j\circ \pi)^{-1} : D_c^b(X,\mathbb C) \to D_c^b(f^{-1}(0), A),
		\end{displaymath}
		where $A := \mathbb C[t_1^{\pm},...,t_r^{\pm}]$, cf. Sabbah \cite{Alexander Module}, see also \cite{BLSW17}.
		Here $D_c^b(X,\mathbb C)$ is the derived category of bounded complexes of sheaves with $\mathbb C$-constructible cohomologies on $X$, and $D_c^b(f^{-1}(0),A)$ is the derived category of bounded complexes of sheaves with $A$-constructible cohomologies on $f^{-1}(0)$. We call $\psi_F\mathbb C_X$ the \textit{specialization complex of $F$}. 	\end{definition}
	
	When $r = 1$, $\psi_f(\mathbb C_X)$ as defined above equals the shift by $-1$ of Deligne's nearby cycles complex together with the action of the monodromy, which is defined using $R\pi_*$ instead of $R\pi_!$, cf. \cite{Bry86}. 
	
	\begin{definition} Let $x \in \{f = 0\}$ and $i\in\bZ$. Let $\Supp^i_x(\psi_F\bC_X)\subset \Spec(A)=(\bC^*)^r$ be the support of $\mathcal H^i(\psi_F\mathbb C_X)_x$ as a finitely generated $A$-module. The {\it local monodromy support}
		is	$\mathrm{Supp}_x(\psi_F\mathbb C_X) := \bigcup_{i\in \mathbb Z}\mathrm{Supp}^i_x(\psi_F\bC_X)$,  a reduced Zariski closed subset of   $(\bC^*)^r$, since this is a finite union.  The \textit{monodromy support locus of $F$} is defined to be $$S(F) := \bigcup_{x \in \{f = 0\}} \mathrm{Supp}_x(\psi_F\mathbb C_X)) \subset (\mathbb C^*)^r.$$  Since we are assuming that condition (\ref{eqCon}) holds, this is a finite union.	
	\end{definition}
	
	When $r=1$, $S(f)$ is the set of monodromy eigenvalues of $f$, that is, numbers which occur as eigenvalues of the monodromy on some cohomology of some Milnor fiber of $f$ at some point in the zero locus of $f$.

	The monodromy support locus admits an interpretation in terms of cohomology support loci of local systems of rank one. Let $U_x$ be the complement of the zero locus of $f$ in a small ball at $x$ in $X$. Consider the restriction of $F$ as a mapping $U_x\to (\bC^*)^r$. For $\lam\in(\bC^*)^r$, let $\cL_\lam$ denote the pullback to $U_x$ of the rank one $\bC$-local system on $(\bC^*)^r$ with monodromy $\lam_i$ around the $i$-th missing coordinate hyperplane. 
	
	\begin{theorem}\label{torsion translated}$\;$
		\begin{enumerate}
			\item (\cite{Bud15}, see also \cite[Prop. 3.7]{LM})
			$\mathrm{Supp}_x(\psi_F\mathbb C_X)=\{\lambda\in(\bC^*)^r\mid H^\ubul(U_x,\cL_\lam)\neq 0\}=\bigcup_i\nu_x^{-1}(\cV^i_1(U_x))$, the last expression using  the notation from Proposition \ref{propRefW}. 
			\item (\cite{BW15})
			$\mathrm{Supp}_x(\psi_F\mathbb C_X)$ is a conjugation-invariant finite union of torsion-translated complex affine subtori of $(\mathbb C^*)^r$.
		\end{enumerate}
	\end{theorem}

	The  specialization complex $\psi_F\bC_X$ admits a good behaviour under specializations of $F$, as we describe next.

	\begin{definition}\label{defNDGSp}
		For $M = (m_{kj}) \in \mathbb N^{p\times r}$, the \textit{specialization} $F^M$ of $F$ is the $p$-tuple of holomorphic functions $(\prod_{j=1}^r f_{j}^{m_{1j}},..., \prod_{j=1}^r f_{j}^{m_{pj}})$. The specialization $F^M$ is called \textit{non-degenerate} if the  induced map on tori $(\mathbb C^*)^r \to (\mathbb C^*)^p$ given by $M$ is surjective and $\sum_{k=1}^p m_{ki} \neq 0$ for all $i$ such that $f_i^{-1}(0)\neq \emptyset$.
	\end{definition}
	\begin{remark}
		If $p = 1$ and $m_{11},...,m_{1r} >0$, the specialization $F^M=\prod_{j=1}^rf_j^{m_{1j}}$ is always non-degenerate.
	\end{remark}
	
	For the  following see  \cite[2.3.8]{Alexander Module}, \cite[Proposition 3.31]{Bud15} :
	
	\begin{proposition}\label{specialization of alexander module}
		If $G = F^M$ is a non-degenerate specialization of $F$, then for all $x \in \{f = 0\} \subset X$, we have
		$			(\tau_M)^{-1}(\mathrm{Supp}_x(\psi_F\mathbb C_X)) = \mathrm{Supp}_x(\psi_G(\mathbb C_X)),
		$		where
		$ \tau_M :  (\mathbb C^*)^p \to (\mathbb C^*)^r$ is given by $(\lambda_1,...,\lambda_p) \mapsto (\prod_{k=1}^p \lambda_k^{m_{k1}},..., \prod_{k=1}^p \lambda_k^{m_{kr}})$.
	\end{proposition}
	
	We will also need  the module-theoretical refinement of this statement, cf. \cite[2.3.10]{Alexander Module}. We point out that the original version by Sabbah is stated only for $x \in \bigcap_{i=1}^r \{f_i = 0\}$, but the proof  works for all $x\in \{f = 0\}$.
	
	\begin{proposition}\label{module-theoretical specialization}
		Let $G = F^M$ be a non-degenerate specialization of $F$. Let $A_r$, $A_p$ be the affine coordinate rings of $(\mathbb C^*)^r$ and $(\mathbb C^*)^p$, respectively. View $A_p$ is as an $A_r$-algebra via the morphism $\tau_M$.	
		Then for all $x \in \{f = 0\} \subset X$, there is a natural quasi-isomorphism
		$			(\psi_G \mathbb C_X)_x \xa{\sim} (\psi_F \mathbb C_X)_x \otimes_{A_r}^L A_p
		$ of complexes of $A_r$-modules.	\end{proposition}

	\subs{\bf Bernstein-Sato ideals and zeta functions.}\label{B-S polynomial} We continue with the same setup as in \ref{Alexander module and its specialization}.
	
	\begin{definition}
		The {\it local  Bernstein-Sato ideal} of $F$ at $x \in \{f = 0\}$ is the ideal $B_{F,x}$ generated by all polynomials $b\in \mathbb C[s_1,...,s_r]$ such that 
		\be\label{eqBS}		b(s_1,...,s_r) \cdot f^{s_1}_1 \cdots f_r^{s_r} \in \mathscr D_{X,x}[s_1,...,s_r] \cdot f_1^{s_1+1} \cdots f_r^{s_r+1},
		\ee	where $\mathscr D_{X,x}\cong \bC\{x_1,\ldots,x_n\}[\pa/\pa x_1,\ldots, \pa/\pa x_n]$ is the stalk of the sheaf of linear differential operators with holomorphic coefficients of $X$ at $x$. Here $s_i$ are independent variables and the equation should be interpreted formally, using the usual rules of differentiation on the right-hand side. 
	\end{definition}
	
	The  Bernstein-Sato ideal is non-zero, by Sabbah \cite{Sab87}. In the case $r = 1$, the monic generator of $B_{f,x}$ is called the local Bernstein-Sato polynomial of $f$ at $x$.
	
	\begin{definition}
		The {\it  Bernstein-Sato ideal of $F$} is $B_F:=\bigcap_{x\in\{f=0\}}B_{F,x}\subset \bC[s_1,\ldots,s_r]$.  The condition (\ref{eqCon}) implies that as $x$ varies, there are only finitely many local Bernstein-Sato ideals $B_{F,x}$ possible by \cite[Cor. 3, Rem. 4]{BMM}. So $B_F$ is a non-zero ideal. We denote by $Z(B_{F})$  the zero locus in $\bC^r$.		\end{definition}

	If $X$ is a smooth affine complex variety and $f_i$ are polynomials, then $B_F$ agrees with the  purely algebraic definition obtained by replacing in (\ref{eqBS}) $\mathscr D_{X,x}$ with the ring of algebraic linear differential operators on $X$.

	The literature on the Bernstein-Sato ideals $B_F$ is too extensive to list it here exhaustively. They were initially also considered in \cite{May, Bio, BrMay, BMM1, BMM2, Bah} in restricted setups, among other papers;  \cite{Bud15} opened new connections and many papers appeared after it. There are other ideals of Bernstein-Sato type defined by other versions of  (\ref{eqBS}), see Section \ref{secObl}.

	For a germ of a holomorphic function $F:(\bC^n,0)\to (\bC^r,0)$, $F$ is smooth if and only if $B_F=((s_1+1)\ldots(s_r+1))$, by \cite[Prop. 1.2]{BrMay}.

	As it is well-known, Bernstein-Sato ideals provide information
	about the polar locus $P(Z_F)$ in $\bC^r$ of the multivariable archimedean zeta function $Z_F$, see \cite[Theorem 5.3.2]{Igu} for the case $r=1$.
	
	\begin{proposition}\label{B-S proof of one direction}
		Let $F = (f_1,...,f_r) : X \to \mathbb C^r$ be an analytic mapping from a complex manifold $X$  satisfying (\ref{eqCon}). Then:
		
		\begin{enumerate}
			\item The polar locus $P(Z_F)$ is a union of $\bZ^r$-translates of hyperplanes in $Z(B_F)$.

			\item		Let  $m\in\bN_{>0}$ and let $C$ be a polar hyperplane of $Z_{F}$ in $\bC^r$. Let $m\in\bN$ be such that $C$ intersects the half-open box $(-m-1,-m]^r$.
			Then $C +k\cdot \bm 1 \subset Z(B_F)$ for some $k\in\{0,\ldots,m\}$,  where $\bm 1 = (1,...,1) \in \mathbb N^r$. 
			
			\item The polar order of $Z_F$ along the $C$ in (2) is at most $\sum_{k = 0}^m \mathrm{ord}_{C+k\cdot \bm 1}(B_F)$, where $\mathrm{ord}$ means the vanishing order of an ideal at a prime divisor. 
			
			\item $\mathrm{Exp}(P(Z_{F})) \subset \mathrm{Exp}(Z(B_F))$ in $(\bC^*)^r$.
			
		\end{enumerate}
	\end{proposition}
	\begin{proof}
		Assume first that  $X = U \subset \mathbb C^n$ is an open subset and $f_1,...,f_r \in \mathcal O_{\mathbb C^n}(U)$. Then  $\omega \in \Omega^{n,n}_{c}(U)$ can be uniquely written as $\omega = \varphi \, \mathrm{d} \bm z$, where $\mathrm{d} \bm z := \mathrm{d} z_1 \wedge \mathrm{d} \bar z_1 \wedge \cdots \wedge \mathrm{d} z_n \wedge \mathrm{d} \bar z_n$ and $\varphi \in C_c^{\infty}(U)$.  By using a partition of unity, one  reduces the problem of  meromorphically extending  $Z_{F,\omega}(s_1,\ldots,s_r)$ beyond $\mathrm{Re}(s_i)> 0$ for all $i$, to a finite collection of open subsets of $U$, denoted by $\{U_k\}$, such that on each $U_k$ there exists $P_k \in \mathscr D_X(U_k)[s_1,...,s_r]$ and $0\neq b_k \in \mathbb C[s_1,...,s_r]$ satisfying $b_k(s_1,...,s_r) \cdot f^{s_1}_1 \cdots f_r^{s_r} = P_k \cdot f_1^{s_1+1} \cdots f_r^{s_r+1}$. We can assume $b_k$ are as in Theorem \ref{Zero_Loci_of_B-S} (2).
		
		Now we drop  $k$ from the notation. We have
		\begin{displaymath}
			b(s_1,\ldots,s_r) \cdot \int_{U} \prod_{j=1}^r\vert f_j\vert^{2s_j} \varphi\, \mathrm{d}\bm{ z} = \int_{U} \prod_{j=1}^r(P f_j^{s_j+1}){\bar f}_j^{s_j}  \varphi\, \mathrm{d}\bm{ z}.
		\end{displaymath}
		The meromorphic extension follows from integration by parts and repeating this procedure as in the proof  for $r=1$ from \cite[Theorem 5.3.2]{Igu}, taking  $p=0$ there.		 We conclude that 
		if $b(s_1,\ldots,b_s)\in B_F$ is a polynomial as in Theorem \ref{Zero_Loci_of_B-S} (2), then $P(Z_F)\subset \bigcup_{k\in\bN}(Z(b)-k\cdot \bm 1)$. The global case under the assumption (\ref{eqCon}) of theis statement reduces to this case by using a partition of unity for the forms $\omega$.
		Since the polar locus $P(Z_F)$ is locally the zero locus of an analytic hypersurface, it follows that each irreducible component $C$ of $P(Z_F)$ is algebraic and equal to a translate of a hyperplane in $Z(B_F)$, and this shows (1). Meanwhile, the polar order is controlled by vanishing orders along translates of $C$ of $B_F$. The amount of translation is limited by the above process, which gives parts (2) and (3). Part (4) follows from  (1). 	\end{proof}
	
	We record an observation about propagation of poles for archimedean zeta functions:
	
	\begin{lemma}\label{propagation_of_poles}
		Let $F$ be as in Proposition \ref{B-S proof of one direction}. If $H$ is a polar hyperplane of $Z_F$ of order $m$, then $H-\bm e_i$ is polar hyperplane of $Z_F$ of order $\geq m$ for all $i=1,\ldots,r$. 
	\end{lemma}
	\begin{proof}
		It suffices to  consider the case $i = 1$. Suppose $H$ is an order $m$ polar hyperplane of $Z_{F,\omega}$, where $\omega \in \Omega_c^{n,n}(X)$. Then $H-\bm e_1$ is an order $m$ polar hyperplane of $Z_{F,\omega}(s_1+1,s_2,...,s_r) = Z_{F,\vert f_1\vert^2 \omega}(s_1,...,s_r)$, and hence is an order $\geq m$ polar hyperplane of $Z_F$.
	\end{proof}

	\section{Cohomological determinantal-factor ideals}\label{Cohomological determinantal-factor ideal}
	
	\subs{\bf Definition and properties.}\label{General definition}
	We now define and draw the main properties of the cohomological determinantal-factor ideals of a complex of modules. We show that they generalize the notion of ``maximal size of Jordan blocks''. This will be applied to monodromy support loci.

	Let $R$ be a Noetherian ring. Let $E^\ubul$ be a bounded above complex of $R$-modules such that its cohomology $H^i(E^\ubul)$ is a finitely generated $R$-module for all $i\in\bZ$. Then $E^\ubul$ admits a bounded-above free resolution $F^\ubul$ of finite rank in each degree, that is, a quasi-isomorphism $ F^\ubul \xa{\sim} E^\ubul$, cf. \cite[III.12.3]{Hartshorne}. Denote by $d^\ubul$ the differential of $F^\ubul$. 	Let $I_m(d^{i-1})$ denote the  ideal generated by the $m\times m$ minors of $d^{i-1} : F^{i-1}\to F^i$. If $m \le 0$, respectively $m > \min\{\mathrm{rank}(F^{i-1},\mathrm{rank}(F^i)\})$, we set $I_m = (1)$, respectively $I_m = 0$.

	\begin{lemma} Fix $i,k\in\bZ$. Set
		$$\mathcal I_k^i(E^\ubul) :=I_{r_i(F^\ubul)-k}(d^{i-1})\subset R,$$  where
		$r_i(F^\ubul) := \sum_{l \geq i} (-1)^{l-i}\mathrm{rank}(F^l)$.  Then $\mathcal I_k^i(E^\ubul)$ does not depend on the choice of $F^\ubul$. 
	\end{lemma}
	\begin{proof} The proof is similar to \cite[Definition-Proposition 2.2]{BW15jump}.
		Let $G^\ubul \xa{\sim} E^\ubul$ be another bounded-above resolution of $E^\ubul$ of finitely generated free modules.  Using the same reduction in loc. cit., it suffices to consider the case when $R$ is local and $G^\ubul = F^\ubul \oplus M$, where $M$ is a shift of the complex $0 \to R \overset{1}{\longrightarrow} R \to 0$. 
		Suppose this complex is added on the $i$-th position, that is,
		\begin{displaymath}
			G^\ubul = \big[\cdots \to F^{i-2} \to F^{i-1} \oplus R \xa{d^{i-1}\oplus 1} F^i\oplus R \to F^{i+1} \to \cdots \big]
		\end{displaymath}
		By  definition, only $\mathcal I^{i-1}_k,\mathcal I^{i}_k$, and $\mathcal I_k^{i+1}$ might change. However, we compute that
		\begin{align*}
			\mathcal I^{i-1}_k(G^\ubul)  = I_{r_{i-1}(G^\ubul)-k}\begin{pmatrix} d^{i-2} \\ 0\end{pmatrix},\,
			\mathcal I^i_k(G^\ubul) & = I_{r_i(G^\ubul)-k}\begin{pmatrix} d^{i-1} & 0 \\ 0 & 1\end{pmatrix},\,
			\mathcal I^{i+1}_k(G^\ubul)  = I_{r_{i+1}(G^\ubul)-k}\begin{pmatrix} d^{i} & 0\end{pmatrix}.
		\end{align*}
		Since $r_{i-1}(G^\ubul) = r_{i-1}(F^\ubul)$, $r_i(G^\ubul) = r_i(F^\ubul)+1$, and $r_{i+1}(G^\ubul) = r_{i+1}(F^\ubul)$, we see that these ideals are the same as for $F^\ubul$. 
	\end{proof}
	
	\begin{definition}
		We call $\mathcal I_k^i(E^\ubul) $ the \textit{cohomological determinantal-factor ideals} of $E^\ubul$.	
	\end{definition}
	
	If $E^\ubul = N[-i]$ is the shift of an $R$-module $N$, then $\mathcal I_k^i(E^\ubul)= I_{n-k}(\varphi)=\mathrm{Fitt}_k(N)$, the  $k$-th Fitting ideal of $N$, where $R^m \xa{\varphi} R^n \to N \to 0$ is a finite presentation of $N$. Hence cohomological determinantal-factor ideals generalize Fitting ideals.
	
	\begin{corollary}
		If $E^\ubul$ is quasi-isomorphic to $E'^\ubul$, then $\mathcal I_k^i(E^\ubul) = \mathcal I_k^i(E'^\ubul)$. That is,  $\mathcal I_k^i$ is well defined on the derived category $D_{fg}^+(R)$ of bounded-above complexes of $R$-modules with finitely generated cohomologies.
	\end{corollary}
	
	\begin{corollary}\label{base change of cohomological determinantal-factor ideal}
		Let $R$ and $E^\ubul$ be  as above, and let $S$ be a Noetherian $R$-algebra. Then $\mathcal I_k^i(E^\ubul) \cdot S = \mathcal I_k^i(E^\ubul \otimes^L_R S)$, where $-\otimes_R^L S$ is the derived tensor product.
	\end{corollary}
	\begin{proof}
		We may assume $E^\ubul$ is a bounded-above complex of free $R$-modules of finite ranks. Both sides are generated by the images of the determinantal ideals of differentials of $E^\ubul$ in $S$, so the equality holds.
	\end{proof}

	When $R$ is a field, by definition $\mathcal I_k^i(E^\ubul) = 0$ if $\sum_{l\geq i} (-1)^{l-i} \dim H^l(E^\ubul) \geq k+1$ and $\mathcal I_k^i(E^\ubul) = (1)$ if $\sum_{l\geq i} (-1)^{l-i} \dim H^l(E^\ubul) \leq k$. Thus:
	
	\begin{corollary}\label{determinantal factors jump loci}
		For any maximal ideal $\mathfrak m$ of $R$, $\mathcal I_k^i \subset \mathfrak m$ if and only if $$\sum_{l\geq i} (-1)^{l-i} \dim_{R/\mathfrak m} H^l(E^\ubul \otimes_R^L R/\mathfrak m) \geq k+1.$$ 
	\end{corollary} 
	Hence a geometric interpretation of the cohomological determinantal-factor ideals is that they place a closed subscheme structure on the locus of jumps of alternating cohomology dimensions as above.

	\begin{corollary}\label{Iki Fitting}
		If $R$ is a principal ideal domain and all cohomologies of $E^\ubul$ are torsion, then $\mathcal I_k^i(E^\ubul) = \mathrm{Fitt}_k(H^i(E^\ubul))$. In particular, $\cI^i_{-1}(E^\ubul)=0$ and $\mathcal I_0^i(E^\ubul)\neq 0$.
	\end{corollary}
	\begin{proof}
		We can assume that $E^\ubul$ is a bounded-above complex of free $R$-modules of finite rank. Let $Z^i := \ker d^i$ and $B^i := \mathrm{im}\, d^{i-1}$. Then $Z^i$ and $B^i$ are  free $R$-modules of finite rank since both are submodules of $E^i$.  Furthermore, $d^i : E^i/Z^i \xa{\simeq} B^{i+1}$ is an isomorphism of free modules. So we have a decomposition $E^i = Z^i \oplus B^{i+1}$. Thus $E^\ubul = \bigoplus_{i\in \mathbb Z} [0 \to B^{i} \xa{d^{i-1}} Z^{i} \to 0]$, where $Z^i$ is added on the $i$-th position, which we use to compute $\mathcal I_k^i(E^\ubul)$.
		We use the free presentation $B^i \xa{d^{i-1}} Z^i \to H^i(E^\ubul) \to 0$ to compute $\mathrm{Fitt}_k(H^i(E^\ubul))$. To show the claim, it suffices to show that $r_i(E^\ubul) = \mathrm{rank}(Z^i)$. Since $H^j(E^\ubul)$ are torsion for all $j\in \mathbb Z$, we have $\mathrm{rank}(Z^j) = \mathrm{rank}(B^j)$. Then we have $r_i(E^\ubul) = \sum_{j \geq i} (-1)^{j-i}\mathrm{rank}(E^j) = \sum_{j\geq i} (-1)^{j-i}\big(\mathrm{rank}(B^{j+1})+\mathrm{rank}(j^i)\big) = \mathrm{rank}(B^i)$.
	\end{proof}
	
	\begin{remark}\label{rmkLinTr}
		Let $\varphi : \mathbb C^m \to \mathbb C^m$ be a linear transformation. Then the morphism of $\bC$-algebras $\mathbb C[t] \to \mathrm{End}(\mathbb C^m)$, $t\mapsto \varphi$, equips $V = \mathbb C^m$ with a $\mathbb C[t]$-module structure with a presentation
		$V\otimes \bC[t]\xa{t\cdot \mathrm{id}_m-\varphi} V\otimes\bC[t]\to V\to 0.	
		$
		Let $b_k$ be the monic generator of the ideal $\cI^0_{k}(V)=I_{m-k}(t\cdot \mathrm{id}_m-\varphi)=\mathrm{Fitt}_{k}(V)\subset \bC[t]$. The polynomials $b_k$ are by definition the classical {\it determinantal factors} of $\varphi$, the successive quotients of which give the non-zero diagonal elements of the Smith normal form of $t\cdot \mathrm{id}_m-\varphi$. 	In particular, $b_{0}/b_{1}$ is the minimal polynomial of $\varphi$. The multiplicities of the factors of $b_{0}/b_{1}$ are the  maximal sizes of Jordan blocks for the corresponding eigenvalues. 
	\end{remark}

	Following this remark we formalize a generalization of the concept of  maximal size of Jordan blocks. We  assume now that $R$ is a regular integral domain. For a prime ideal $\mathfrak p\subset R$ 
	of height one and a non-zero ideal $I\subset R$, let $v_{\mathfrak p}(I)$ denote the order of $I_{\mathfrak p}$ in the discrete valuation ring $R_{\mathfrak p}$, and let  
	$$\mathrm{div}(I) = \sum_{\mathrm{ht}\, \mathfrak p = 1} v_{\mathfrak p}(I) \cdot V(\mathfrak p) \in \mathrm{Div}(\mathrm{Spec}\, R)
	$$
	be the divisor associated to $I$.
	
	\begin{definition}\label{determinantal divisor} Let $R$ be a regular integral domain.
		Assume that $E^\ubul \in D_{fg}^+(R)$ has torsion cohomologies. For $i,k\in\bZ$ with $k\geq 0$:
		\begin{itemize}
			\item The \textit{determinantal-factor divisors} of $E^\ubul$ are 
			$\Delta_{k}^i(E^\ubul) := \mathrm{div}(\mathcal I_{k}^i(E^\ubul)) \in \mathrm{Div}(\mathrm{Spec}\, R). $
			\item
			The \textit{minimal divisors} of $E^\ubul$ are
			$			M^i(E^\ubul) := \Delta_{0}^i(E^\ubul)-\Delta_{1}^i(E^\ubul).
			$		
		\end{itemize}
	\end{definition}
	
	\begin{remark}
		Under the torsion assumption each ideal $\mathcal I_{k}^i(E^\ubul)$ with  $k \geq 0$, is non-zero. To see this we can localize the complex at an arbitrary height one prime ideal $\mathfrak p$, then $\mathcal I_k^i(E^\ubul)_{\mathfrak p} = \mathcal I_k^i(E_{\mathfrak p}^\ubul) \neq 0$ by Corollary \ref{Iki Fitting}. Hence $\mathcal I_k^i(E^\ubul) \neq 0$ and the above definitions make sense.
	\end{remark}
	
	The following proposition relates $\mathcal I_k^i$ and the annihilator of the $i$-th cohomology.
	\begin{proposition}\label{specialization of CIDF of H^i}
		Let $R$ be a regular integral domain.
		Assume that $E^\ubul \in D_{fg}^+(R)$ has torsion cohomologies. Then $M^i(E^\ubul) = \mathrm{div}(\mathrm{Ann}_R(H^i(E^\ubul)))$, in particular $M^i(E^\ubul)$ is an effective non-zero divisor. Moreover,
		$(\Delta_0^i(E^\ubul))_{\mathrm{red}} = (\mathrm{div}(\mathrm{Ann}_R(H^i(E^\ubul))))_{\mathrm{red}}$.
	\end{proposition}
	\begin{proof}
		Since cohomologies, annihilators, and cohomological determinantal-factor ideals are all compatible with localization, we can replace $R$ with $R_{\mathfrak p}$, where $\mathfrak p$ is height one prime ideal. The rest follows from Corollary \ref{Iki Fitting} and Lemma \ref{PID Fitting and Annihilator}.
	\end{proof}
	\begin{lemma}\label{PID Fitting and Annihilator}
		Suppose $R$ is a PID, $N$ is a finitely generated torsion $R$-module, and $\mathrm{Fitt}_j(N) = (b_j)$. Then $\mathrm{Ann}_R(N) = (b_0/b_1)$ and $V(\mathrm{Ann}_R(N)) = V(b_0)$.
	\end{lemma}
	\begin{proof}
		We can decompose $N$ as $N = \bigoplus_{j=0}^m R/(a_j)$, with $0\neq (a_0) \subset (a_1) \subset ... \subset (a_m)$. Then $R^{m+1} \xa{(a_0,...,a_m)} R^{m+1} \to N\to 0$ is a presentation. We find that $b_j = \prod_{k\geq j} a_k$ and $\mathrm{Ann}_R(N) = (a_0)$. We thus conclude the proof.
	\end{proof}

	\subs{\bf Comparison with cohomology jump ideals.}\label{Comparison between determinantal and jump}
	Let $R$ be a Noetherian ring and $E^\ubul$ be a bounded-above complex of $R$-modules with finitely generated cohomologies.  
	The \textit{cohomology jump ideals} of $E^\ubul$ are defined by 
	\begin{displaymath}
		\cJ_k^i(E^\ubul) := I_{\mathrm{rank}(F^i)-k+1}(d^{i-1}\oplus d^i)\subset R
	\end{displaymath}
	where $(F^\ubul, d^\ubul)$ is a bounded-above complex of finite-rank free modules quasi-isomorphic to $E^\ubul$, see \cite{BW15jump}. The ideals
	$\cJ_k^i(E^\ubul)$ are also independent of the choice of $F^\ubul$, are compatible with base change,  generalize the Fitting ideals of a finitely generated module, and
	place a closed subscheme structure on the set of maximal ideals $\mathfrak m\subset R$ with $\dim H^i(E^\ubul\otimes^L_RR/\mathfrak m)\ge k$. 
	The cohomological determinantal-factor ideals  refine  the cohomology jump ideals:	
	\begin{proposition} We have
		$			\cJ_j^i(E^\ubul) = \sum_{k+l = j-1} \mathcal I_{k}^i(E^\ubul)\cdot \mathcal I_{l}^{i+1}(E^\ubul).
		$		Consequently, for every height one prime ideal $\mathfrak p\subset R$, $
		v_{\mathfrak p}(\cJ_j^i(E^\ubul)) = \min_{k+l = j-1}\{v_{\mathfrak p}(\mathcal I_{k}^i(E^\ubul)) + v_{\mathfrak p}(\mathcal I_{l}^{i+1}(E^\ubul))\}.
		$	\end{proposition}
	\begin{proof}
		Replacing $E^\ubul$ by a free resolution $F^\ubul$ as above, we have from definitions that $\cJ_j^i(E^\ubul) = \sum_{(k,l) \in T_{j,i}} \mathcal I_{k}^i(E^\ubul)\cdot \mathcal I_{l}^{i+1}(E^\ubul)$, where $T_{j,i} = \{(k,l) \in \mathbb Z^2 \mid (r_i-k) + (r_{i+1}-l) = \mathrm{rank}(F^i)-j+1\}$. Since $r_i + r_{i+1} = \mathrm{rank}(F^i)$, we conclude the proof.
	\end{proof}

	\section{Local systems and the specialization complex}\label{secLocos}
	In this section we  apply the above concepts  to place a scheme structure on  monodromy support loci. We prove here Proposition  \ref{rigidity theorem of cohomological determinantal-factor jump loci}, Proposition \ref{propRefW}, and Theorem \ref{zero locus of cohomological determinantal-factor ideals}.
	
	\subs We first prove Proposition  \ref{rigidity theorem of cohomological determinantal-factor jump loci}. Let $Y$ be a topological space with homotopy type of a finite CW complex. Recall, for example from \cite{BW15jump}, that the moduli space $\mathcal M_B(Y,m)$ is the GIT quotient of the representation space $\mathcal R_B(Y,m)=\mathrm{Hom}(\pi_1(Y,y),\mathrm{GL}_m(\mathbb C))$ by the conjugation action of $\mathrm{GL}_m$, where $y\in Y$ is a fixed point.	Since $\pi_1(Y,y)$ is finitely presented, $\mathcal R_B(Y,m)$ and $\mathcal M_B(Y,m)$ are affine schemes of finite type over $\bC$. The closed points of $\mathcal M_B(Y,m)$ are the isomorphism classes of semi-simple $\bC-$local systems of rank $m$ on $Y$. A representation corresponding to a closed point in $\mathcal R_B(Y,m)$ maps in $\mathcal M_B(Y,m)$ to its semi-simplification. One defines the cohomology jump loci set-theoretically as
	\begin{displaymath}
		\tilde{\mathcal V}_k^i := \{\rho \in \mathcal R_B(Y,m) \mid \dim_{\mathbb C} H^i(Y, L_\rho) \geq k\},
	\end{displaymath}
	where $ L_\rho$ is the local system corresponding to the representation $\rho$. Then ${\mathcal V}_k^i :=\cV^i_k(Y,m)$ is the image of $\tilde{\mathcal V}_k^i $ in $\mathcal M_B(Y,m)$.

	A closed subscheme structure is placed on these sets as follows, cf. \cite{BW15jump}.
	Let $R$ be the coordinate ring of $\mathcal R_B(Y,m)$ and let $\mathcal L$ be the rank $m$ $R$-local system corresponding to the universal rank $m$ representation of $\pi_1(Y,y)$. That is, $\mathcal L$ is a local system of free rank $m$ $R$-modules on $Y$ such that $\mathcal L_{\rho} = \mathcal L\otimes_{R} R/\mathfrak m_{\rho}$ for all $\rho \in \mathcal R_B(Y,m)$ with maximal ideal $\mathfrak m_\rho \subset R$.	Let $a : Y \to \{pt\}$ be the map to  a point, then $Ra_*(\mathcal L)$ is represented by a bounded complex
	of  $R$-modules with finitely generated cohomologies. For $\rho \in \mathcal R_B(Y,m)$, by the projection formula, we have $Ra_*(\mathcal L) \otimes^L_R R/\mathfrak m_\rho = Ra_*(\mathcal L\otimes_RR/\mathfrak m_\rho) = Ra_*(\mathcal L_\rho)$, which computes the cohomologies of $\mathcal L_\rho$. Hence the cohomology jump ideal $\cJ^i_k(Ra_*(\mathcal L)) \subset R$ places a closed subscheme structure on $\tilde{\mathcal V}_k^i $.  By general properties of the GIT quotient map, the image $\mathcal V^i_k$ of $\tilde {\mathcal V}^i_k$ is also closed, and thus $\cJ^i_k(Ra_*(\mathcal L))$ induces a closed subscheme structure on $\mathcal V^i_k$.

	\begin{proof}[Proof of Proposition \ref{rigidity theorem of cohomological determinantal-factor jump loci}]		Similarly, by Corollary \ref{determinantal factors jump loci} we see that
		the cohomological determinantal-factor ideal $\cI^i_{k-1}(Ra_*(\mathcal L))$ places a closed subscheme structure on $\mathcal A^i_k:={\mathcal A}_k^i(Y,m)$, by first defining it in the representation space $\mathcal R_B(Y,m)$. This proves (1). 
		
		Set-theoretically, $\mathcal A^i_k$ is in the boolean algebra generated by the cohomology jump loci, that is, it is obtained from  $\{\mathcal V^j_l\}_{j,l}$ by finite unions, intersections, and taking complements. On the other hand, this boolean algebra is generated by torsion-translated complex affine subtori of $\mathcal M_B(Y,1)$ in the case $m=1$, if $Y$ is a smooth $\bC$-algebraic variety \cite{BW20}, or the small ball complement of a  $
		\bC$-analytic subset \cite{BW15}. Since $\mathcal A^i_k$ are Zariski closed, this shows (2).
	\end{proof}

	\begin{remark} A generalization and a different proof of Proposition \ref{rigidity theorem of cohomological determinantal-factor jump loci} are as follows.	
		Since cohomological determinantal-factor ideals are compatible with base change,  \cite[Theorem 5.14.3 (3)]{BW20} holds with $h^i$ replaced by the function $r_i:D^b_c(pt,\bC)\to \bZ$, $r_i(E^\ubul)=\sum_{l\ge i}(-1)^{l-i}\dim H^i(E^\ubul)$.  As a consequence, \cite[Lemma 6.2.2]{BW20} holds for $r_i$, that is, with the terminology from there, $r_i$ is an absolute $\bQ$-semi-continuous function. Then Theorems 6.4.3 (4),  7.3.2 (3), 9.1.1 (3) of \cite{BW20} hold with $h^i$ replaced by $r_i$  as well. So any composition $\phi$ of functors as in \cite[Theorem 6.4.3 (1)]{BW20}  (resp. as in \cite[Theorem 6.4.3 (3)]{BW20}) on smooth complex algebraic varieties  starting from $\mathcal M_B(Y,m)$ and ending with the function $r_i$ in $\bZ$, the inverse image $\phi^{-1}(\bZ_{\ge k})$ is an absolute $\bQ$-constructible (resp. absolute $\bQ$-closed) subset of local systems. In particular, if the rank $m=1$, such a set is in the boolean algebra generated by the (resp. is a finite union of) torsion-translated complex affine subtori of $\mathcal M_B(Y,1)$, by \cite[Theorem 1.3.1]{BW20}.
	\end{remark}

	\subs		Let now $F = (f_1,...,f_r) : X \to \mathbb C^r$ be an analytic mapping from a complex manifold $X$ of dimension $n$ satisfying (\ref{eqCon}). Consider the specialization complex $\psi_F\bC_X$ in $D_c^b(f^{-1}(0), A)$ defined in \ref{Alexander module and its specialization}, where $A = \mathbb C[t_1^{\pm},...,t_r^{\pm}]$.

	\begin{proof}[Proof of Proposition \ref{propRefW}] We recall a fact  due to Sabbah, using the diagram and the notation from \ref{Alexander module and its specialization}.
		Let $\mathcal L =  \tilde \pi_!\mathbb C_{(\mathbb C^*)^r}$. This is the local system of free rank one $A$-modules on $(\bC^*)^r$
		corresponding to the universal rank one representation  of $\pi_1((\mathbb C^*)^r)$, that is, to the homomorphism $\pi_1((\mathbb C^*)^r)\to A^\times$ sending the class of the loop around the $i$-th missing hyperplane to $t_i$. Let $\mathcal L^F$ be the  local system  on $U$ obtained by the inverse image of $\cL$ under $F$. Then by \cite[2.2.8]{Alexander Module}, \cite[Lemma 3.3]{Bud15}, we have $(\psi_F\mathbb C_X)_x \simeq (Rj_*(\mathcal L^F))_x$ as complexes of $A$-modules. 
		
		The zero locus of the ideal $\cJ^i_k$, respectively $\cI^i_{k}$, of this complex parametrizes the set of maximal ideals $\mathfrak m\subset A$ such that the dimension of $H^i((\psi_F\mathbb C_X)_x\otimes^L_A A/\mathfrak m)$ is $\ge k$, respectively the alternative sum of dimensions of cohomologies in degrees $\ge i$ is $\ge k+1$, by Corollary \ref{determinantal factors jump loci}. By the projection formula, $(Rj_*\cL^F)_x\otimes^L_AA/\mathfrak m = (Rj_*(\cL^F\otimes_A^LA/\mathfrak m))_x= (Rj_*(\cL^F\otimes_A A/\mathfrak m))_x = (Rj_*\cL_\lam)_x$, where $\cL_\lam$ is the rank one $\bC$-local system on $U$ obtained from pullback via $F$ from the rank one local system on $(\bC^*)^r$ with monodromy $\lam_i$ around the $i$-th missing hyperplane. Hence $H^i((\psi_F\mathbb C_X)_x\otimes^L_A A/\mathfrak m) = H^i(U_x,\cL_\lam)$. Then parts (1) and (2) of the proposition follow. Part (3) is a consequence of Proposition \ref{propAmu}.
	\end{proof}

	For a point $x \in \{f = 0\}$, the local monodromy support
	$\mathrm{Supp}_x(\psi_F \mathbb C_X)$ is a closed subset of $(\mathbb C^*)^r$ of codimension $>0$, by Theorem \ref{Zero_Loci_of_B-S}. Hence all cohomologies of $(\psi_F\mathbb C_X)_x$ are  torsion $A$-modules and Definition \ref{determinantal divisor} applies: 
	
	\begin{definition}\label{Application to Sabbah}
		For $F$ as above, $x \in \{f = 0\}$, and $i,k\in\bZ$ with $k\ge 0$: 
		\begin{itemize}
			\item the {\it cohomology determinantal-factor ideal} is $\mathcal I_{k}^i(F,x):=\mathcal I_{k}^i((\psi_F\mathbb C_X)_x)\subset A$;		
			\item the {\it determinantal-factor divisor} is  $\Delta^i_{k}(F,x):=\Delta_k^i((\psi_F\mathbb C_X)_x)$; 
			\item the {\it minimal divisor} is $M^i(F,x):=M^i((\psi_F\mathbb C_X)_x)=\Delta^i_0(F,x)-\Delta^i_1(F,x)$;
			\item $\mathrm{ord}^i_C(F,x)$ is the order of a prime divisor $C\subset (\bC^*)^r$ in $M^i(F,x)$.
		\end{itemize}
	\end{definition}

	By   Proposition \ref{propRefW}	 (1) and (2), and Proposition \ref{specialization of CIDF of H^i}, we have:
	
	\begin{prop}\label{propAmu} There is an equality $M^i(F,x)=\divi(\mathrm{Ann}_A(H^i((\psi_F\bC_X)_x)))$
		of non-zero effective divisors in $(\bC^*)^r$. There is an equality of the reduced divisors underlying $\Delta^i_{0}(F,x)$, $M^i(F,x)$,   the codimension-one part of $\Supp^i_x(\psi_F\bC_X)$, and the codimension-one part of $\nu_x^{-1}(\cA^i_{1}(U_x))$.
	\end{prop}
	
	\begin{remark}
		Proposition \ref{propAmu} is further complemented by a relation between $\Supp^i_x(\psi_F\bC_X)$ and the local cohomology Alexander  module of $F$
		by a result of Liu-Maxim \cite[Proposition 3.7]{LM}, which in turn is a correction of \cite[Proposition 3.9]{Bud15}.
	\end{remark}	
	
	\begin{proof}[Proof of Proposition \ref{rmkAmuCon}.] In the case $r=1$ of Proposition \ref{propAmu}, we have $F=f$, $A=\bC[t^\pm]$,   $U_x$ is the complement of $\{f=0\}$ in a Milnor ball at $x$, and $\nu_x:\bC^*\to \mathcal M_B(U_x)$ sends $\lam\in\bC^*$ to the rank one local system $\cL_\lam$ on $U_x$ with monodromy $\lam$ around $\{f=0\}$. Moreover $H^i((\psi_F\bC_X)_x))=H^{i-1}(F_{f,x},\bC)$ where $F_{f,x}$ is the Milnor fiber of $f$ at $x$. So by Proposition \ref{propAmu} there is an equality in $\bC^*$ between the set
		of eigenvalues of monodromy on $H^{i-1}(F_{f,x},\bC)$
		and the set of $\lam$ such that $h^i(U_x,\cL_\lam)-h^{i+1}(U_x,\cL_\lam)+\ldots \ge 1\}.
		$ 
	\end{proof}
	
	\begin{remark} Regarding the equality of sets from Proposition \ref{rmkAmuCon},
		for $i=0$, the left-hand side is trivially empty, and the right-hand side is empty since, $U_x$ being a smooth locally trivial fibration with fiber $F_{f,x}$ over a punctured disc, it has $\chi(U_x)=0$. For $i=1$, and assuming that $x$ is an isolated singularity of $\{f=0\}$, the left-hand side is $\{1\}$, and the right-hand side is $\{\lam\in\bC^*\mid h^0(U_x,\cL_\lam)\ge 1\}$ and therefore also equal to $\{1\}$. In general, a different proof of Proposition \ref{rmkAmuCon} can be given using the Leray spectral sequence for the Milnor fibration.
	\end{remark}

	\begin{proof}[Proof of Theorem \ref{zero locus of cohomological determinantal-factor ideals}]
		Part (1) follows from Proposition \ref{rigidity theorem of cohomological determinantal-factor jump loci} (2) and Proposition \ref{propRefW} (2).
		
		For part (2), let $C = \{h=0\}$ for $h \in A$  irreducible, non-invertible. By part (1) we can assume that $h=\prod_{i=1}^r t_i^{m_i} - \xi$, where $m_1,...,m_r \in \mathbb N$ are coprime, not all zero, and $\xi$ is a unit root. Then there exists a hypersurface given by $g\in A$ such that: (i) $C \not \subset \{g = 0\}$; (ii) all other irreducible components of $S(F)$ are contained in $\{g = 0\}$; (iii) $(\mathcal I_{0}^i(F,x))_g = (h^{c_0})_g$ and $(\mathcal I_1^i(F,x))_g = (h^{c_1})_g$ in the localization $A_g$, and $c_0 > c_1$. It follows immediately that $\mathrm{ord}_{C}^i(F,x) = c_0-c_1$. We show that $U := \{g \neq 0\} \cap C$ has the desired properties.

		
		For a torsion point $(\lambda^{b_1},...,\lambda^{b_r}) \in U$, where $\lambda$ is a unit root and $b_1,...,b_r \in \mathbb N_{>0}$, let $A \to  B:=\bC[t^{\pm}]$ send $t_i$ to $t^{b_i}$ for all $i$. 
		This corresponds to a morphism of varieties $\bC^*\to (\bC^*)^r$ sending $\lam$ to $(\lam^{b_1},\ldots,\lam^{b_r})$.
		Then the image of  $h$ is $t^{\sum_{i=1}^r b_im_i}-\xi$, a polynomial of which
		$\lambda$ is a multiplicity-one root by our assumption. 
		We have equalities of ideals of $B$,
		$$
		\cI^i_k((\psi_F\bC_X)_x)\cdot B\stackrel{(\ref{base change of cohomological determinantal-factor ideal})}{=}
		\cI^i_k((\psi_F\bC_X)_x\otimes^L_AB) \stackrel{(\ref{Iki Fitting})}{=}
		\mathrm{Fitt}_k H^i((\psi_F\bC_X)_x\otimes^L_AB)  \stackrel{(\ref{module-theoretical specialization})}{=}
		\mathrm{Fitt}_k H^i((\psi_{\tilde f}\bC_X)_x)
		$$
		where $\tilde f=\prod_{i=1}^rf_i^{b_i}$. This is further equal to $\mathrm{Fitt}_k H^{i-1}(F_{\tilde f,x},\bC)$, where $F_{\tilde f,x}$ is the Milnor fiber of $\tilde f$ at $x$. Then $c_0-c_1$ is the difference between the  multiplicities of $\lam$ as a root  of the monic generators of $\mathrm{Fitt}_0$ and $\mathrm{Fitt}_1$ of $H^{i-1}(F_{\tilde f,x},\bC)$. The claim  follows by  Remark \ref{rmkLinTr}, or Lemma \ref{PID Fitting and Annihilator}.
	\end{proof}

	\section{Archimedean zeta functions}\label{pf1} 
	
	In this section we prove the rest of the  results from the introduction, Theorems \ref{Thm_Support_Pole},  \ref{Thm_Order controlling},  \ref{Specialization of Archimedean Zeta Function},  as well as the stronger  Theorem \ref{thmGenbc}.

	\subs{\bf Specialization.} 
	Let $X$ be a complex manifold   and $F = (f_1,...,f_r) : X \to \mathbb C^r$ be an analytic mapping satisfying (\ref{eqCon}). Let $\bm s=(s_1,\ldots,s_r)$. We consider a new parameter $s$.	For $\omega \in \Omega_c^{n,n}(X)$  denote by $Z_{F,\omega}(\bm s,s)$  the meromorphic extension to $\bC^{r+1}$  of 
	$$
	(\bm s,s) \mapsto\int_{X} \vert f_1\vert^{2s_1s} \cdots \vert f_r\vert^{2s_rs} \omega,
	$$
	cf. \cite{SZ25}.
	Thus $Z_{F,\omega}(\bs,s)=Z_{F,\omega}(\bs)\circ\rho$, the pullback under the algebraic morphism $\rho:\bC^{r+1}\to\bC^r$, $\rho(\bs, s)= \bs\cdot s$, of the meromorphic function $Z_{F,\omega}(\bs)$. We drop $\omega$ from the notation when we view $Z_{F}(\bs,s)$ as a distribution. A description of the polar locus of $Z_{F}(\bm s,s)$ follows from Proposition \ref{B-S proof of one direction}. Recall  that   the polar locus $P(Z_F)$ is at most a countable union of hyperplanes.
	
	\begin{lemma}\label{lemCrC} The polar locus of $Z_{F}(\bs,s)$ consists of the preimages in $\bC^{r+1}$ of the polar hyperplanes of $Z_F(\bs)$. The polar order of a polar hyperplane $C$ of $Z_{F}(\bs)$ equals the polar order of the polar hypersurface $\rho^{-1}(C)$ of $P(Z_F(\bs,s))$.
	\end{lemma}
	\begin{proof} Clearly, $P(Z_F(\bs,s))\subset \rho^{-1}(P(Z_F(\bs)))$. Conversely, locally, the polar locus of a meromorphic function $\phi$ is the common zero locus of all holomorphic functions which appear as denominators in the representation of $\phi$ as a ratio of two holomorphic functions. Since every component of $P(Z_F(\bs))$ is of the type $C=\{\sum_{i=1}^rc_is_i+c_0=0\}$ with all $c_i\in\bN$ and $c_0>0$, $C$ and $\rho^{-1}(C)$ are locally everywhere irreducible. It follows that $\rho^{-1}(C)$ is a polar hypersurface for $Z_F(\bs,s)$. This proves the first assertion. We have also obtained the polar order of $\rho^{-1}(C)$ in $Z_F(\bs,s)$ is at most the polar order of $C$ in $Z_F(\bs)$. Now,  the polar orders for $Z_F(\bs,s)$ can only decrease under the specialization $s\mapsto 1$, and $P(Z_{F}(\bs,s))\cap\{s=1\}\subset P(Z_F(\bs,1))$. However, $Z_F(\bs,1)=Z_F(\bs)$. This proves the second assertion.
	\end{proof}
	
	For $\bm b\in\bN_{>0}^r$, consider the complex line $\Delta_{\bm b}=\{\bm b\cdot s\mid s\in\bC\}\subset\bC^{r}$ through the origin and $\bm b$. 
	
	\begin{lemma}\label{lemCsa}
		Suppose $C=\{\sum_{i=1}^rc_is_i+c_0=0\}$ is a polar hyperplane of $Z_F(\bs,s)$ of order $m$.  Let $\bm b \in \mathbb N_{>0}^r$ be such that the point $C\cap \Delta_{\bm b}$ in $\bC^r$ lies away from the other polar hyperplanes of $Z_F(\bs)$ and the zero locus of $Z_F(\bs)$. Then $-c_0/\sum_{i=1}^rb_ic_i$ is a pole of order $m$ of $Z_{f_1^{b_1}\cdots f_r^{b_r}}(s)$.\end{lemma}
	\begin{proof} Note that $\Delta_{\bm b}$ is the projection under $\rho$ of the line $L_{\bm b}=\bigcap_{i=1}^r\{s_i=b_i\}\subset\bC^{r+1}$. 
		By Lemma \ref{lemCrC}, the assumption implies that $L_{\bm b}$ intersects the polar hypersurface $\rho^{-1}(C)$ of $Z_{F}(\bs,s)$ away from the other polar hypersurfaces and the zero locus of $Z_F(\bs,s)$. We also have that the polar order of $\rho^{-1}(C)$ is $m$.
		Hence the restriction of $Z_F(\bs,s)$ to $L_{\bm b}$, which equals $Z_{f_1^{b_1}\cdots f_r^{b_r}}(s)$, has $L_{\bm b}\cap\rho^{-1}(C)$ as a pole of order $m$. Parametrizing $L_{\bm b}$ by the coordinate $s$ on $\bC^{r+1}$, we see that $L_{\bm b}\cap\rho^{-1}(C)=\{-c_0/\sum_{i=1}^rb_ic_i\}$.
	\end{proof}

	\begin{proof}[Proof of Theorem \ref{Specialization of Archimedean Zeta Function}.]
		In the  Lemma \ref{lemCsa} we removed a locus of tuples $\bm b$ defined by polar and zero loci of a meromorphic function on $\bC^r$. We are left with an open dense locus in the analytic topology. Moreover, the assumption in the lemma depends only on the directions of the lines $\Delta_{\bm b}$, hence we can rephrase it for $\bm b$ such that $\sum_{i=1}^rb_i=1$. This gives (1).
		
		For (2), we parametrize again $L_{\bm b}\subset \bC^{r+1}$ by the $s$ coordinate. Then  
		the assumption is equivalent to: $\al\in L_{\bm b}$ is a pole of order $m$ of $Z_{F}(\bs,s)_{\mid L_{\bm b}} = Z_{f_1^{b_1}\cdots f_r^{b_r}}(s)$. Since by restriction the pole orders cannot increase, we have that $m$ is smaller than or equal to the sum of the polar orders of the polar hypersurfaces of $Z_{F}(\bs,s)$ containing $\al$. Then the conclusion follows by Lemma \ref{lemCrC}. 
	\end{proof}

	Using  \cite[Theorem 1.7]{SZ25} and Theorem \ref{Specialization of Archimedean Zeta Function} (2),  we get another proof  of the multivariable $n/d$-conjecture of \cite{Bud15} for the case of complete splittings of hyperplane arrangements due to \cite{Wu22}; for another proof see  \cite[Example 4.9]{SV}.

	\begin{theorem} \label{thmNDConj}
		Let $F=(f_1,\ldots f_r)$ be a collection of linear forms in $\bC[x_1,\ldots,x_n]$ such that the associated hyperplane arrangement $f=\prod_{i=1}^rf_i$ is central, essential, and indecomposable, not necessarily reduced.
		Then $\{s_1+\ldots+s_r+n=0\}$ is a polar hyperplane of $Z_F$ and is also contained in the zero locus of $B_F$.
	\end{theorem} 
	\begin{proof}
		Let $H=\{s_1+\ldots+s_r+n=0\}$. Choose a generic point  $\sigma$ of $H\cap \bQ_{<0}^r$. We can assume that $\sigma$ does not lie in any polar hyperplane of $Z_F$ different than $H$, even though at the moment we do not know $H$ is polar hyperplane. Thus we can represent $\sigma$ as $-(n/\sum_ib_i)\cdot\bm b$ with $\bm b=(b_1,\ldots, b_r)\in\bN_{>0}^r$. By \cite[Theorem 1.7]{SZ25}, $-(n/\sum_ib_i)$ is a pole of $Z_{f_1^{b_1}\cdots f_r^{b_r}}(s)$, since $\bm b/\sum_{i=1}b_i$ was chosen generically in $\{s_1+\ldots+s_r=1\}\cap\bQ_{>0}^r$. By Theorem \ref{Specialization of Archimedean Zeta Function} (2), the sum of the polar orders of the polar hyperplanes of $Z_F$ containing $-(n/\sum_ib_i)\cdot\bm b$ is $\ge 1$. By the way we chose $\sigma$, the only possible contribution to this sum is from $H$, and hence $H$ is a polar hyperplane of $Z_F$. The assumptions on the hyperplane arrangement $f$ imply that $r>n$. Hence  by Proposition  \ref{B-S proof of one direction}  there cannot be any other translate of a hyperplane in $Z(B_F)$ that equals $H$. Hence $H\subset Z(B_F)$. 
	\end{proof}

	\subs{\bf Proof of Theorem \ref{Thm_Support_Pole}.}
	By Proposition \ref{B-S proof of one direction} and Theorem \ref{Zero_Loci_of_B-S}, 	$\mathrm{Exp}(P(Z_{F})) \subset S(F)$ and both are finite unions of torsion-translated complex affine subtori of codimension one in $(\bC^*)^r$. We focus on the converse.
	It suffices to show all torsion points of $S(F)$ are contained in $\mathrm{Exp}(P(Z_{F}))$. Pick a torsion point of $S(F)$ and write it as $(\lambda^{b_1},...,\lambda^{b_r})$, where $\lambda = e^{-2\pi i /b}$ and $b,b_1,...,b_r \in \mathbb N$ with $b\neq 0$. By Proposition \ref{specialization of alexander module}, we have that $\lambda$ is a monodromy eigenvalue of $f_1^{b_1}\cdots f_r^{b_r}$. Hence, by Theorem \ref{Triangle of three local invariant} (2), there exists a pole $\alpha$ of $Z_{f_1^{b_1}\cdots f_r^{b_r}}(s)$ such that $e^{2\pi  i\alpha} = \lambda$. By Theorem \ref{Specialization of Archimedean Zeta Function} (2), there exists a polar hyperplane $ C=\{ \sum_{i=1}^r c_i s_i + c_0 = 0\}$ of $Z_{F}(\bm \bs)$, with $c_i\in\bN$ and $c_0\neq 0$, such that $\sum_{i=1}^r c_i b_i \alpha + c_0 = 0$. Thus $(\lambda^{b_1},...,\lambda^{b_r}) = \mathrm{Exp}(b_1\alpha,...,b_r\alpha)\in \Exp(C)\subset\mathrm{Exp}(P(Z_{F}))$. $\Box$

	\subs{\bf Proof of Theorem \ref{Thm_Order controlling}.} We prove a more general statement than Theorem \ref{Thm_Order controlling}:
	
	\begin{theorem}\label{thmGenbc} Let $X$ be a complex manifold of dimension $n$  and $F = (f_1,...,f_r) : X \to \mathbb C^r$ be an analytic mapping satisfying (\ref{eqCon}).
		Let $C$ be an irreducible component of codimension one of $\Supp^i_x(\psi_F\bC_X)$. 
		\begin{enumerate}
			
			\item Let $(e^{2\pi i\al b_1},\ldots,e^{2\pi i\al b_r})$ with $\alpha \in (-1,0]\cap\bQ$ and $\bm b=(b_1,\ldots,b_r) \in \mathbb N_{>0}^r$, be a point of $U$ as in Theorem \ref{specialization of multiplicitiy}(2). Then there exists $q\in\bN$ and a polar hyperplane $H$ of $Z_F$ of polar order $\ge \mathrm{ord}_{C}^i(F,x)$ such that $C=\Exp(H)$ and $\bm b\cdot(\al-q)\in H$.

			\item If in addition $i\ge 2$, one can take $q=i-1$.
			
			\item If in addition $i\ge 2$ and $\al=0$, then there exists $q\in\bN$ and a polar hyperplane $H$ of $Z_F$ of polar order $\ge \mathrm{ord}_{C}^i(F,x)+1$ such that $C=\Exp(H)$ and $-q\cdot\bm b\in H$.
			
		\end{enumerate}
	\end{theorem}	
	\begin{proof}
		We write our point of $U$ as $\sigma = (\lambda^{b_1},\ldots,\lambda^{b_1})$, where $\lambda = e^{2\pi i\alpha}$ is the unit root. Let $\tilde f =  f_1^{b_1}\ldots f_r^{b_r}$. 
		By   Theorem \ref{specialization of multiplicitiy} (2), the maximal size of Jordan blocks with eigenvalue $\lam$ of the monodromy on the $(i-1)$-th cohomology of the Milnor fiber of $\tilde f$ at $x$ equals $\mathrm{ord}_{C}^i(F,x)$. Note that $\mathrm{ord}_{C}^i(F,x)\ge 1$ by	Proposition \ref{propAmu}.

		Theorem \ref{Triangle of three local invariant} (3) (a) implies that there exists $q\in\bN$ such that $\al-q$ is a pole of $Z_{\tilde f}(s)$ of 	 order $\geq \mathrm{ord}_{C}^i(F,x)$. 	If $i\ge 2$	 then Theorem \ref{Triangle of three local invariant} (3) (b) implies that one can take $q=i-1$, that is,
		$\al-(i-1)$ is a pole of $Z_{\tilde f}(s)$ of order $\geq \mathrm{ord}_{C}^i(F,x)$.
		
		In any case,
		by Theorem \ref{Specialization of Archimedean Zeta Function}	(2), the sum of the polar orders of the polar hyperplanes of $Z_F$ containing  $\bm b\cdot(\al-q)$ is $\ge \mathrm{ord}_{C}^i(F,x)$. 
		
		Note that $\Exp(\bm b\cdot(\al-q))=\sigma$. Let $H=\{\sum_{i=1}^rc_is_i+c_0=0\}$ be a lift  of $C$ under $\Exp$  containing 	$\bm b\cdot(\al-q)$, where $c_i\in\bN$  and $c_0\neq 0$. By Theorem \ref{Thm_Support_Pole}, $C\subset \Exp(P(Z_F))$. However no  $\bZ^r$-translate of a polar hyperplane of $Z_F$ different than $H$ can contain the lift  $\bm b\cdot(\al-q)$ of $\sigma$, by the way we chose $\sigma$. Since $\mathrm{ord}_{C}^i(F,x)\ge 1$, $H$ must be a polar hyperplane of $Z_F$ containing  $\bm b\cdot(\al-q)$, it is the unique such one, and moreover its polar order in $Z_F$ is  $\ge \mathrm{ord}_{C}^i(F,x)$. This proves  (1) and (2).
		
		For part (3), Theorem \ref{Triangle of three local invariant} (3) (c) implies that there exists $q\in\bN$ such that $-q$ is a pole of $Z_{\tilde f}(s)$ of  order $\ge \mathrm{ord}_{C}^i(F,x)+1$. The above argument then gives a polar hyperplane $H$ of $Z_F$ containing $-q\cdot\bm b$	of polar order $\ge \mathrm{ord}_{C}^i(F,x)+1$.
	\end{proof}

	\begin{remark}$\;$
		\begin{itemize}		
			\item Theorem \ref{thmGenbc} (1) (respectively (2), (3)) is a generalization of Theorem \ref{Triangle of three local invariant} (3) (a) (respectively (b), (c)).
			
			\item Even in the case $r=1$ these present a small generalization, by allowing $\bm b =b$ to be different than 1, compared with Barlet's results. For example, Theorem \ref{thmGenbc} (2) for $r=1$ says: 
			
			{\it If $j\ge 1$ and $e^{2\pi i\al b}$ with $\al\in (-1,0]\cap \bQ$ and $b\in \bN_{>0}$ is a monodromy eigenvalue on the cohomology $H^{j}$ of the Milnor fiber of $f$ at $x$ with a Jordan block of size $m\ge 1$, then  $\al b-jb$ is a pole of $Z_f(s)$ of  order $\ge m$.}
			
			\item In the last remark we see that, compared with  Theorem \ref{Triangle of three local invariant} (3) (c), we  guarantee the existence of infinitely many poles of order $\ge m$, not only one pole. This propagation of polar hyperplanes holds also for $r>1$ according to Theorem \ref{thmGenbc} or Lemma \ref{propagation_of_poles}.
		\end{itemize}
	\end{remark}

	\section{Oblique slopes}\label{secObl}

	Now we consider all  ideals of Bernstein-Sato type from \cite[\S 4]{Bud15}. In this section we prove Theorem \ref{thmOBLIQUE} and provide a few complementary results on these more general Bernstein-Sato ideals.
	
	Again, let $F = (f_1,...,f_r) : X \to \mathbb C^r$ be an analytic mapping from a complex manifold $X$ of dimension $n$ satisfying (\ref{eqCon}).
	
	\begin{definition}\label{defBFM} Let $p\in\bN_{>0}$ and $M=(m_{kj})\in\bN^{p\times r}$. The {\it Bernstein-Sato ideal associated to $F$ and $M$ at $x\in\{f=0\}$} is the ideal $B_{F,x}^M$ generated by all $b\in\bC[s_1,\ldots,s_r]$ such that
		\be\label{eqBSM}		b(s_1,...,s_r) \cdot f^{s_1}_1 \cdots f_r^{s_r} \in \sum_{k=1}^p \mathscr D_{X,x}[s_1,...,s_r] \cdot f_1^{s_1+m_{k1}} \cdots f_r^{s_r+m_{kr}}.
		\ee
		We set $B_F^M=\bigcap_{x\in\{f=0\}}B_{F,x}^M$. 
		If $M=I_r\in\bN^{r\times r}$ is the identity matrix, we set  $B_F^I=B_F^{I_r}$.  
		
	\end{definition}

	Note that $B_F=B_F^{\,\bm{1}}$ where $\bm{1}=(1,\ldots,1)\in \bN^r$. If $M$ contains a zero row, or if $M$ contains a non-zero row $\bm m_k$ such that, for all $1\le j \le r$, $f_j$ is invertible if $m_{kj}\neq 0$, then $B_F^M=(1)$. In all other cases, $B_F^M$ contains a non-constant polynomial in $\bC[s_1,\ldots,s_r]$.	The ideals $B_F^I$ were  considered  in \cite{May, UC} and briefly in \cite{Bud15}, but they are largely unexplored.

	We recall the following results.
	Let $t_j$ be the ring isomorphism of $\bC[s_1,\ldots, s_r]$ defined by $t_j(s_i) = s_i +\delta_{ij}$, for $j=1,\ldots, r$. Denote by $T_j$ the associated isomorphism  $\bC^n$ defined by $T_j(\bm\alpha)=\bm\al-\bm e_j$ for a vector $\bm\al\in\bC^n$.
	
	\begin{theorem}\label{thmBFMs} Let $\bm m\in\bN^r$. Then:
		\begin{enumerate}
			\item (\cite[\S 4]{Bud15})  There are inclusions of ideals of $\bC[s_1,\ldots,s_r]$:
			$$
			\prod_{\substack{1\le j\le r\\ m_j>0}} \prod_{k=0}^{m_j-1} t_1^{m_1}\ldots t_{j-1}^{m_{j-1}}t_j^k \,B_F^{\bm e _j}\subset B_F^{\bm m} \subset \bigcap_{\substack{1\le j\le r\\ m_j>0}}\bigcap_{k=0}^{m_j-1} t_1^{m_1}\ldots t_{j-1}^{m_{j-1}}t_j^k \, B_F^{\bm e_j}.
			$$
			Both sides can be changed by a permutation $\pi$ of $\{1,\ldots, r\}$, see \cite[Remark 4.10]{Bud15}.
			In particular, for all $\pi$ there is an equality of zero loci in $\bC^r$:
			$$
			Z(B_F^{\bm m})=\bigcup_{\substack{1\le j\le r\\ m_{\pi(j)}>0}}\bigcup_{k=0}^{m_{\pi(j)}-1} T_{\pi(1)}^{m_{\pi(1)}}\ldots T_{\pi(j-1)}^{m_{\pi(j-1)}}T_{\pi(j)}^k \,Z(B_F^{\bm e_{\pi(j)}}).
			$$

			\item (\cite{BS2}) Every irreducible component of $Z(B_F^{\bm m})$ of codimension 1 is a  hyperplane of type $\{c_1s_1+\ldots+c_rs_r+c_0=0\}$ with $c_i\in\bN$, $c_0>0$, and for each such hyperplane there exists $1\le i\le r$ with $m_i\ne 0$ such that $c_i>0$. 
			\item (\cite{BS2}) Every irreducible component of  $Z(B_F^{\bm m})$ of codimension $>1$  can be translated by an element of $\bZ^r$ inside a component of codimension 1.
			
			\item (\cite{BS2}) $$\Exp (Z(B_F^{\bm m}))=\cup_{j}\Exp (Z(B_F^{\bm e_j})) =\cup_{j} S_j(F)$$ where the unions are over $1\le j\le r$ with $m_j\ne 0$ and $f_j$ not invertible, and
			$S_j(F):=\Supp_{(\bC^*)^r}(\psi_F\bC_X)|_{f_j^{-1}(0)}$.

		\end{enumerate}
	\end{theorem}
	
	Parts (2)-(4) generalize Theorem \ref{Zero_Loci_of_B-S} (2)-(4) by taking $\bm m=\bm 1$, and one has $S(F)=\cup_j S_j(F)$. Part (4) gives the precise topological information contained by $Z(B_F^{\bm m})$ for all vectors $\bm m\in\bN^r$.

	\begin{lemma}\label{lemPrep}
		Let $M\in \bN^{p\times r}$ be a matrix with row vectors $\bm m_k\neq \bm 0$, $1\le k\le p$. Then: 
		
		\begin{enumerate}
			
			\item  There are inclusions of ideals $\sum_{k=1}^p B^{\bm m_k}_{F}\subset B_{F}^{\,M}\subset B_F^I $. In particular, $$Z(B_F^I)\subset Z(B_{F}^{\,M})\subset\cap_{k=1}^p Z(B^{\bm m_k}_{F}).$$
			
			\item There is an inclusion of ideals $\prod_{k=0}^{\max(M)-1}(t_1\ldots t_r)^k B_F \subset B_F^M$, where $\max(M)$ is the maximum of the entries of $M$. In particular,
			$$
			Z(B_F^M)\subset \cup_{k=0}^{\max(M)-1} (Z(B_F)-k\cdot\bm 1).
			$$
			
			\item Let $a\in\bN_{>0}$ and let $M_a$ be a matrix whose rows consist of all $\bm m\in\bN^r$ such that $|\bm m|:=\sum_{j=1}^rm_j=a$. Then
			$$ \prod_{\bm m'\in\bN^r \, :\, \vert \bm m'\vert \le a-1} t_1^{m'_1}\ldots t_r^{m'_r}B_F^{I} \subset B_F^{M_a}.$$
			In particular, 
			$Z(B_F^{M_a})\subset Z(B_F^I)-\bN^r$.
			
			\item (\cite{May}) For $1\le j\le r$, let $F_j=(f_1,\ldots,\hat f_j,\ldots, f_r)$. Denote by $\check B_{F_j}\subset \bC[s_1,\ldots,s_r]$ the ideal generated by $B_{F_j}\subset\bC[s_1,\ldots,\hat s_j,\ldots, s_r]$. Then there exists $a\in\bN$ such that
			$$
			(\bigcap_{j=1}^r \check B_{F_j})\cdot\prod_{\bm m\in\bN^r\,:\, |\bm m|\le a}t_1^{m_1}\ldots t_r^{m_r}B_F^I\subset B_F.
			$$
			In particular, 
			$$
			\bigcup_{j=1}^r Z(\check B_{F_j})\cup \bigcup_{\bm m\in\bN^r\,:\, |\bm m|\le a} (Z(B_F^I)-\bm m)\supset Z(B_F).
			$$
		\end{enumerate}
		
	\end{lemma}
	\begin{proof}
		(1) follows directly from the definition.  (2) is obtained by noticing that $b_1(\bm s)b_2(\bm s+\bm 1)...b_l(\bm s+l\cdot \bm 1) \bm f^{\bm s} \in \mathscr D_{X}[\bm s] \cdot \bm f^{\bm s+(l+1)\cdot\bm 1}  \subset \sum_{k=1}^p \mathscr D_{X}[\bm s]\cdot \bm f^{\bm s+\bm m_{k}}$ for all $b_1,\ldots ,b_l \in B_F$ and $l \ge \max(M)-1$. From the inclusions of ideals in (1) and (2), only $B_{F}^{\,M}\subset B_F^I $ uses that $M$ has no zero rows.
		For (3) one applies inductively (\ref{eqBSM}) by increasing  $|\bm m'|$. Part (4) is essentially \cite[Proposition 1.1]{May}.
	\end{proof}

	\begin{remark} 
		We see that $B_F^I$ is the largest of all  ideals of Bernstein-Sato type $B_F^M$ with $M$ having no zero rows, equivalently, the zero locus $Z(B_F^I)$ is the smallest.	
		The ideal $B_F^I$ is unreasonably easy to compute even in cases when $B_F$ is not computable, see Example \ref{exaBM2}. We do not have an explanation for this. We note that $Z(B_F^I)$ could have codimension $>1$, unlike $Z(B_F)$, see \cite[Example 4.20]{Bud15}. In fact, $Z(B_F^I)=\emptyset$ if and only if any of the $f_j$ is invertible.
	\end{remark}

	\begin{lemma}\label{lemObBFI} The codimension one part of the zero locus $Z(B_F^M)$ is a union of hyperplanes of type  $\{c_1s_1+\ldots+c_rs_r+c_0=0\}$ with $c_i\in\bN$, $c_0>0$. All hyperplanes in $Z(B_F^I)$ are oblique. \end{lemma}
	\begin{proof}
		This follows from Lemma \ref{lemPrep} (1) and Theorem \ref{thmBFMs} (2).
	\end{proof}
	
	\begin{lemma}\label{applem2}
		Let $H$ be an oblique polar hyperplane of $Z_F$. Then $H \subset Z(B_F^{I}) - \mathbb N^r$. 
	\end{lemma}
	\begin{proof} It is enough to prove it for $X=\bC^n$.
		By our assumption, there exists $\varphi\in C_c^{\infty}(\mathbb C^n)$ such that $H$ is a polar hyperplane of $Z_{F,\varphi}$. Suppose $H$ is not contained in $ Z(B_F^{I}) - \mathbb N^r$. By Lemma \ref{lemPrep} (3), for each $a \in \mathbb N_{>0}$, there exists $b \in B_F^{M_a}$ such that $H \not\subset \{b=0\}$. As in the proof of Proposition \ref{B-S proof of one direction}, 
		\begin{equation}
			b(s_1,...,s_r)\int_{\mathbb C^r} \vert f_1\vert^{2s_1}...\vert f_r\vert^{2s_r} \varphi \, \mathrm{d} \bm z = \sum_{\vert \bm m\vert = a} \int_{\mathbb C^n} f_1^{s_1+m_1} ... f_r^{s_r+m_r} \bar f_1^{s_r} ... \bar f_r^{s_r} \phi_{\bm m}\, \mathrm{d} \bm z, 
			\label{GBSM}
		\end{equation}
		where $\phi_{\bm m} \in C_c^{\infty}(\mathbb C^n)[s_1,...,s_r]$. We observe that $\int_{\mathbb C^n} f_1^{s_1+m_1} ... f_r^{s_r+m_r} \bar f_1^{s_r} ... \bar f_r^{s_r} \phi_{\bm m}\, \mathrm{d} \bm z$ is holomorphic on the region $U_{\bm m} := \{\bm s \in \mathbb C^r \mid \mathrm{Re}(s_i) + m_i/2 > 0\text{ for all }i\}$. Take $a \gg 0$ such that each $U_{\bm m}$ intersects with $H$. Note that the right-hand side of (\ref{GBSM}) is a $\mathbb C[s_1,...,s_r]$-combination of archimedean zeta functions, so (\ref{GBSM}) is also an equation for meromorphic functions. However, $H$ is not a polar hyperplane for each term of the right-hand side, so it is not a polar hyperplane of $Z_{F,\varphi}$, a contradiction. 
	\end{proof}

	\begin{proof}[Proof of Theorem \ref{thmOBLIQUE}.]  (1) This is proved in Lemma \ref{lemObBFI}.

		(2) By Lemma \ref{lemPrep} (1)-(2) we have $Z(B_F^{I}) \subset Z(B_F^M) \subset Z(B_F)-\mathbb N\cdot {\bm 1}$. Looking at the oblique parts only and exponentiating, we get	
		\be\label{eqIncl}\Exp(P(Z_F)_{ob})\subset\Exp(Z(B_F^{I})_{ob}) \subset \Exp(Z(B_F^M)_{ob}) \subset \Exp(Z(B_F)_{ob})=S(F)_{ob}\ee
		where the first inclusion follows from Lemma \ref{applem2}, and  the last inclusion follows from Theorem \ref{Zero_Loci_of_B-S} (4).	But $\mathrm{Exp}(P(Z_F)_{ob}) = S(F)_{ob}$ by Theorem \ref{Thm_Support_Pole}. Hence all inclusions in (\ref{eqIncl}) are equalities. 
		
		(3) We also have now that a non-oblique hyperplane $H\subset Z(B_F)$ must 
		come from the zero loci $Z(B_{F_j})$, by Lemma \ref{lemPrep} (4). With the same argument, now we apply the last lemma  inductively.
	\end{proof}

	We complement this section by recording proving the following on specializations and Bernstein-Sato type ideals:
	
	\begin{theorem} (\cite[Conjecture 4.25]{Bud15})
		Let $M\in\bN^{p\times r}$ define a non-degenerate specialization of $F$, see Definition \ref{defNDGSp}, which we denote by $F^M=(\prod_{j=1}^r f_{j}^{m_{1j}},..., \prod_{j=1}^r f_{j}^{m_{pj}})$. Let $\bm m\in \bN^p$. Let $ \tau_M :  (\mathbb C^*)^p \to (\mathbb C^*)^r$ be given by $(\lambda_1,...,\lambda_p) \mapsto (\prod_{k=1}^p \lambda_k^{m_{k1}},..., \prod_{k=1}^p \lambda_k^{m_{kr}})$. Then
		$$
		{\tau}_M^{-1} (\Exp (Z(B_F^{\,\bm m\cdot M}))) = \Exp (Z(B_{F^M}^{\bm m})).
		$$
		In particular, for all $M$ with nonzero columns, 
		$
		\tau^{-1}_M (\Exp(Z(B_F))) = \Exp (Z(B_{F^M})).
		$
	\end{theorem}
	\begin{proof} Denote by $\bm c_j$ the columns of $M$, for $1\le j\le r$.
		Using Theorem \ref{thmBFMs} (4), we  have 
		$$
		{\tau}_M^{-1} (\Exp (Z(B_F^{\,\bm m\cdot M}))) =   \tau_M^{-1} \big(\bigcup_{j\, :\, \bm m\cdot \bm c_{j}\neq 0}\Supp_{(\bC^*)^r}(\psi_F\bC_X|_{Z(f_j)})\big) =$$
		$$= \tau_M^{-1} (\Supp_{(\bC^*)^r}(\psi_F\bC_X|_{\cup_{j\, :\, \bm m\cdot \bm c_{j}\neq 0}Z(f_j)})),
		$$
		$$
		\Exp (Z(B_{F^M}^{\bm m})) = \bigcup_{k\,:\, \bm m_k\neq 0} \Supp_{(\bC^*)^p}(\psi_{F^M}\bC_X|_{Z(F^M_k)}) = \Supp_{(\bC^*)^p}(\psi_{F^M}\bC_X|_{\cup_{k\,:\, \bm m_k\neq 0}Z(F^M_k)}).
		$$
		Since $\cup_{j\, :\, \bm m\cdot \bm c_{j}\neq 0}Z(f_j) = \cup_{k\,:\, \bm m_k\neq 0}Z(F^M_k)$, the first claim follows from Proposition \ref{specialization of alexander module} for non-degenerate specializations. The second claim follows from the first by taking $\bm m=\bm 1$. Of course, the second claim also follows directly from Proposition \ref{specialization of alexander module} and Theorem \ref{Zero_Loci_of_B-S} (4).
	\end{proof}
	
	\begin{remark}\label{rmkConjs}$\;$
		\begin{enumerate}
			\item It is conjectured in \cite{Bud15} that the ideal $B_F$ is generated by products of polynomials of  type $c_1s_1+\ldots+c_rs_r+c_0$ with $c_i\in\bN$, $c_0>0$. We conjecture the same for the ideals $B_F^M$ for all matrices $M$.  The subconjecture that the codimension $>1$ components of the zero loci $Z(B_F)$, and more generally $Z(B_F^M)$, are linear is still open. 
			\item  We still lack a topological interpretation of the whole $\Exp(Z(B_F^M))$ for $M$ with more than one row, including the case $M=I_r$.
			\item We conjecture that $B_F$, and more generally $B_F^{\bm m}$ for $\bm m\in\bN^r$, can be recovered from the ideals $B_F^{\bm e_j}$ using any permutation $\pi$ of $\{1,\ldots, r\}$ as follows:
			$$B_F^{\bm m} \stackrel{?}{=} 
			\bigcap_{\substack{1\le j\le r\\ m_{\pi(j)}>0}}\bigcap_{k=0}^{m_{\pi(j)}-1} t_{\pi(1)}^{m_{\pi(1)}}\ldots t_{\pi(j-1)}^{m_{\pi(j-1)}}t_{\pi(j)}^k \, B_F^{\bm e_{\pi(j)}}.
			$$
			The inclusion $\subset$, and  the equality of their radicals, hold by Theorem \ref{thmBFMs} (1). It was already remarked  in \cite{Bud15} that   no counterexamples were found.
		\end{enumerate}
	\end{remark}

	\section{Examples}\label{secExa}

	\begin{example}\label{exaBM} (\cite[5.1]{BM}) Let $F:\bC^3\to \bC^2$ be $F=(x^2z-y^3,z)$. The Bernstein-Sato ideal is 
		$$\textstyle B_F=((\prod_{k=5,6,7}(6s_1+k))(s_2+1)(\prod_{k=4,5,7,8}(3s_1+3s_2+k))),$$
		using \cite{Lev, Oa}. Therefore
		$$\textstyle
		S(F)=\Exp(P(Z_F))=Z((\prod_{k=1,5,6}(t_1-e^{2\pi ik/6 }))(t_2-1)(\prod_{k=1,2}(t_1t_2-e^{2\pi ik/3})))$$
		as a subset of $(\bC^*)^2$ with coordinates $t_1$, $t_2$, by Theorem \ref{Thm_Support_Pole}.
		Hence the  slopes of the polar locus of $Z_F$ are $(1,0),(0,1),(1,1)$. This confirms that there is an oblique slope for the polar locus of $Z_F$, as found by Barlet-Maire \cite[5.1]{BM}. However, this corrects the slope found in \cite[5.1]{BM} and it shows that $(1,1)$ is the unique oblique slope. 	
		
		We also have:
		$$B_F^I=((3s_1 + 3s_2 + 4)(3s_1 + 3s_2 + 5))\cap  (s_1 + 1, s_2 + 1),$$
		$$B_F^{\bm e_1}=((s_1+1)(6s_1+5)(6s_1+7)(3s_1+3s_2+4)(3s_1+3s_2+5)),$$
		$$ B_F^{\bm e_2} =(s_2+1,(3s_1+3s_2+4)(3s_1+3s_2+5)).$$
		This illustrates the various results in this section, including Theorem \ref{thmOBLIQUE}.
		
		Note that one cannot apply \cite[Corollary 1.10]{Bud15} to produce an oblique slope here, since the topological Euler characteristic of the projective complement of $\{f=0\}$ is 0.
	\end{example}

	\begin{example}\label{exaBM2} (\cite[5.2]{BM}) For $F=(x^4+y^4+x^2yz,z)$, there is an oblique slope $(4,1)$ for $B_F$, since $\chi(\mathbb P^2\setminus Z(f))=4\neq 0$, by using \cite[Corollary 1.10]{Bud15}. Hence  this is also an oblique slope for the polar locus of $Z_F$, by  Theorem \ref{Slopes of zeta function}, which confirms \cite[5.2]{BM}. In fact, computing a log resolution we  see that $(4,1)$ is the only oblique slope.  For the computation of the Euler characteristic we used Aluffi's algorithm implemented in \cite{Al, Ch}. 
		
		Alternatively, we can use  Bernstein-Sato type ideals. For this example, $B_F$ could not be calculated by a computer. However:
		$$\textstyle B_F^I=(4s_1 + 5, s_2 + 1)\cap (2s_1 + 3, s_2 + 1)\cap (s_1 + 1, (s_2  +  1)^2)\cap ( \prod_{k=3}^5(4s_1 + s_2 +k)).$$
		It follows that $$S(F)_{ob}=\Exp(P(Z_F)_{ob})=Z(t_1^4t_2-1)\subset (\bC^*)^2$$ by Theorem \ref{thmOBLIQUE}. This implies again that $(4,1)$ is the unique oblique polar slope, and in fact it gives more information about the location of the oblique components.
		
		In fact we can do better using Theorem \ref{thmBFMs} (1). The computer gives:
				$$\textstyle B_F^{\bm e_1} =((s_1+1)^2)(\prod_{k=5,7,9,11}(8s_1+k))(\prod_{k=3}^8(4s_1+s_2+k))),$$
		$$\textstyle B_F^{\bm e_2} =((s_2+1)\prod_{k=3}^5(4s_1+s_2+k)).$$
		These ideals determine the zero locus of $B_F$ by Theorem \ref{thmBFMs} (1). We obtain in $\bC^2$:
		$$\textstyle
		Z(B_F) =Z((s_1+1)(s_2+1)(\prod_{k=5,7,9,11}(8s_1 + k)) (\prod_{k=3}^9( 4s_1 + s_2 +k))),
		$$
		although we could not compute $B_F$. By  Theorem \ref{Thm_Support_Pole}, 
		we obtain in  $(\bC^*)^2$:
		$$\textstyle
		\Exp(P(Z_F))=S(F)= \Exp(Z(B_F)) = Z((t_1-1)(t_2-1)(t_1^4+1)(t_1^4t_2-1)). 
		$$
		According to Remark \ref{rmkConjs} (3), $B_F$ is conjecturally equal to
		$$B_F^{\bm e_1}\cap t_1B_F^{\bm e_2}=\textstyle\big((s_2 + 1)(s_1  + 1)^2 \cdot \big(\prod_{k=5,7,9,11}(8s_1+k) \big)
		\cdot \big(\prod_{k=3}^9 (4s_1+s_2+k)\big)\big).$$
	\end{example}

	\begin{example} (\cite{Bah}) Let $F=(x_1^a+x_2^b,x_1^c+x_2^d)$ with $a,b,c,d\in\bN_{>0}$, $b,c>1$, and $bc>ad$.
		If $a=d=1$, then $B_F=((s_1+1)(s_2+1))$ by \cite[Thm. 1]{Bah}. Thus in this case $\Exp(P(Z_F))=Z((t_1-1)(t_2-1))$ in $(\bC^*)^2$. If $a\ge 2$ or $d\ge 2$, then the   slopes of $\pi_2(W_{F,0}^\#)$ as in Theorem \ref{Slopes of zeta function} (6) are $(1,0)$, $(0,1)$, $(b,d)$, $(a,c)$ by \cite[Prop. 2]{Bah}. Hence these are also the slopes of $P(Z_F)$  by Theorem \ref{Slopes of zeta function}.
	\end{example}

\end{document}